\documentclass[preprint]{elsarticle}

\usepackage[english]{babel}
\usepackage{graphicx}% Include figure files
\usepackage{dcolumn}% Align table columns on decimal point
\usepackage{bm}% bold mathManuscript Title:\\

\usepackage{upgreek}
\usepackage{amsmath}
\usepackage{amssymb} %for mathbb{G}
\usepackage{diagbox}
\usepackage{commath}
\usepackage{natbib}
\usepackage{amsmath}
\usepackage{algorithm}
\usepackage[noend]{algpseudocode}

\newcommand*\rfrac[2]{{}^{#1}\!/_{#2}}

\usepackage{upgreek}
\usepackage{amsmath}
\usepackage{amssymb} %for mathbb{G}
\usepackage{bbm} %for matrix 1
\renewcommand{\t}{\mathbf}
\usepackage{multirow}

\usepackage{diagbox}
\usepackage{commath}

\usepackage[utf8]{inputenc}
\usepackage{xcolor}

\usepackage[hidelinks]{hyperref} %non colored ref links

\renewcommand{\t}{\mathbf}
\usepackage{mathtools}

\usepackage{layouts} %to measure text width in inches

\usepackage[aboveskip=1pt, belowskip=1pt]{subcaption}
\usepackage{pifont}% http://ctan.org/pkg/pifont
\newcommand{\cmark}{\ding{51}}%
\newcommand{\xmark}{\ding{55}}%

\usepackage{glossaries}
\newacronym[longplural={representative volume elements}]{rve}{RVE}{representative volume element}
\newacronym{fem}{FEM}{finite element method}
\newacronym{fe}{FE}{finite element}
\newacronym{xfem}{XFEM}{extended finite element method}
\newacronym{fg}{FG}{Fourier-Galerkin}
\newacronym{ldpm}{LDPM}{lattice discrete particle model}
\newacronym{fft}{FFT}{fast Fourier transform}
\newacronym{srm}{SRM}{spectral \gls{rve} method}
\newacronym{gpu}{GPU}{graphics processing unit}
\newacronym{cpu}{CPU}{general purpose processor}
\newacronym{ufl}{UFL}{unified form language}
\newacronym{pde}{PDE}{partial differential equation}
\newacronym{dry}{DRY}{\emph{Don't repeat yourself!}}
\newacronym{hpc}{HPC}{high-performance computing}
\newacronym{scitas}{SCITAS}{Scientific IT and application support}
\newacronym{asr}{ASR}{alkali-silica reaction}
\newacronym{aar}{AAR}{alkali-aggregate reaction}
\newacronym{pi}{PI}{principal investigator}
\newacronym{api}{API}{application programming interface}
\newacronym{fefe}{FE$^2$}{FEM squared}
\newacronym{ofen}{SFOE}{Swiss Federal Office of Energy}
\newacronym{empa}{EMPA}{Swiss Federal Laboratories for Materials Science and Technology}
\newacronym{lmc}{LMC}{Laboratory of Construction Materials}
\newacronym{lsms}{LSMS}{Computational Solid Mechanics Laboratory}
\newacronym{epfl}{EPFL}{Ecole Polytechnique Fédérale de Lausanne}
\newacronym{cadd}{CADD}{coupled atomistics and discrete dislocations}
\newacronym{lammm}{LAMMM}{Laboratory for Multiscale Mechanics Modeling}
\newacronym{kit}{KIT}{Karlsruher Institut für Technologie}
\newacronym{cms}{IAM--CMS}{Institute for Applied Materials -- Computational Materials Science}
\newacronym{rom}{ROM}{reduced order model}
\newacronym{cg}{CG}{conjugate gradient}
\newacronym{sla}{SLA}{sequential linear algorithm}
\newacronym{itz}{ITZ}{Interfacial transition zone}
\newacronym{csh}{CSH}{Calcium silicate hydrate}
\newacronym{psd}{PSD}{positive semi definite}
\newacronym{spsd}{SPSD}{symmetric positive semi definite}
\newacronym{pnsd}{NPSD}{symmetric negative semi definite}
\newacronym{snpsd}{SNPSD}{symmetric non-positive semi definite}
\newacronym{nsd}{NSD}{negative semi definite}
\newacronym{bfgs}{BFGS}{Broydon-Fletcher-Goldfarb-Shanno (algorithm)}
\newacronym{lbfgs}{LBFGS}{Low memory Broydon-Fletcher-Goldfarb-Shanno}
\newacronym{dof}{DoF}{degree of freedom}
\newacronym{faief}{FAIEF}{first order approximation incremental energy functional}

\newcommand{\HTranspose}[1]{#1^\mathrm{*}}
\newcommand{\test}{\vek{\zeta}}
\newcommand{\higherorders}{\mathcal{O}(||\vek{p}_i||^3)}

\newcommand{\fluc}[1]{\widetilde{#1}}
\newcommand{\compat}[1]{\Check{#1}}
\newcommand{\Fouri}[1]{\widehat{#1}}

\usepackage{boites, boites_exemples}
\newcounter{till_comment_counter}

\newcommand{\pucc}{\mathrm{\Omega}_0}

\newcommand{\vek}[1]{\mathchoice{\displaystyle\boldsymbol{#1}}
  {\textstyle\boldsymbol{#1}}{\scriptstyle\boldsymbol{#1}}
  {\scriptscriptstyle\boldsymbol{#1}}}

\newcommand{\stress}{\ensuremath{\vek{\sigma}}}
\newcommand{\tangent}{\ensuremath{\mathbb{B}}}

\newcommand{\strain}{\ensuremath{\vek{\varepsilon}}}
\newcommand{\nablastrain}{\nabla_{\strain}}

\newcommand{\FourierTrans}{\mathcal{F}}
\newcommand{\iFourierTrans}{\mathcal{F}^{-1}}

\newcommand{\Doper}{\mathcal{D}}

\newcommand*{\rom}[1]{\expandafter\romannumeral #1}

\begin{document}
% title needs revision
\title{Non-convex, ringing-free, FFT-accelerated solver using an incremental
  approximate energy functional }

\address[1]{Department of Mechanical Engineering,
  \'Ecole Polytechnique F\'ed\'erale de Lausanne, 1015 Lausanne, Switzerland}

\address[2]{Department of Microsystems Engineering, University of Freiburg,
  Georges-K\"ohler-Allee 103, 79110 Freiburg, Germany}

\address[3]{Faculty of Civil Engineering, Czech Technical University in Prague,
  Th\'akurova 7, 166 29 Prague 6, Czech Republic}

\author[1]{Ali Falsafi\corref{cor1}}%{A. Falsafi}
\ead{ali.falsafi@epfl.ch}
\author[2]{Richard J. Leute}% {R. Leute}
\author[3]{Martin Ladecký}%{M. Ladecký}
% \author[3]{Jan Zeman}%{J. Zeman}
% \author[3]{Ivana Pultarová}%{I. Pultarová}
% \author[2]{Lars Pastewka}%{L. Pastewka}
\author[1]{Till Junge}%{T. Junge}
\cortext[cor1]{Corresponding author}
\date{Jun 2022}

\begin{abstract}
  Fourier-accelerated micromechanical homogenization has been developed and
  applied to a variety of problems, despite being prone to ringing
  artifacts. In addition, the majority of Fourier-accelerated solvers applied
  to~\gls{fft}-accelerated schemes only apply to convex problems. We here introduce a
  \gls{faief} that allows to employ modern efficient and non-convex iterative
  solvers, such as trust-region solvers or \gls{lbfgs} in a
  \gls{fft}-accelerated scheme. These solvers need the explicit energy functional of
  the system in their standard form. We develop a modified trust region solver,
  capable of handling non-convex micromechanical homogenization problems such
  as continuum damage employing the~\gls{faief}. We use the
  developed solver as the solver of a ringing-free~\gls{fft}-accelerated
  solution scheme, namely the projection based scheme with finite element
  discretization.
\end{abstract}


\maketitle
\glsresetall
%%%%%%%%%%%%%%%%%%%%%%%%%%%%%%%%%%%%%%%%%%%%%%%%%%%%%%%%%%%%%%%%%%%%% 
%                                                                   % 
% Introduction                                                      %
%                                                                   % 
%%%%%%%%%%%%%%%%%%%%%%%%%%%%%%%%%%%%%%%%%%%%%%%%%%%%%%%%%%%%%%%%%%%%% 
\section{Introduction}
Mechanical homogenization, motivated by the idea of representing a
heterogeneous micro-structure as an equivalent homogeneous medium, aims to
calculate effective mechanical properties of micro-structures, including
homogenized elastic constants and the stress-strain response given the
micro-structure and constitutive laws of the individual components.
For simple micro-structures (e.g. micro-structures containing only linear
elastic phases) the effective properties of a heterogeneous
material can be estimated analytically~\cite{Budiansky1965, Mori1973,
  Norris1985, Hill1985,  Milton2003, Milton1995, Nemat-Nasser2013}.
However, when the micro-structure of a material gets
more complex, analytical methods are generally no longer suitable for the
determination of the effective properties.

Computational homogenization, on the other hand, is an effective
method in up-scaling the behavior of complex micro-structures specially those
consisting of,
\rom{1}. highly nonlinear, or
\rom{2}. evolving phases~\cite{Hill1963, Geers2010}.
Computational homogenization methods are
based on the construction of a micro-scale boundary value problem, the so-called
cell problem, discretizing the solution domain and
solving the governing equation, equilibrium equation for instance,
using numerical schemes such as~\gls{fem}~\cite{Schroeder2014} or
spectral methods~\cite{Moulinec1994, Eisenlohr2013}.

One of the primary applications of computational homogenization is in multi-scale
simulations, where it enables resolving the full micro-structure and studying
the influence of parameters at microscale on the solution of a structural
problem~\cite{Schroeder2014, Liu2016, Schneider2016, Matous2017}.
In a multiscale scheme, such as \gls{fefe}, discretization points (quadrature
points) at the macro-scale are each represented by a \gls{rve} in which the
micro structure of the underlying phases is captured. The strains at
macro-scale discretization points are imposed as boundary conditions on the
micro-scale model and the resultant mean stress and effective tangent extracted
from the solution of the cell problem is passed to the macro-scale. In a
multi-scale approach, usually, the computational costs of the cell problem
solution outweighs that of the higher scale problem since in a single
macroscale load step the~\gls{rve} solver is called at least once per material
point. Accordingly, it is of high importance to optimize and speed-up the
\gls{rve} solution. In comparison with~\gls{fem}, the computational efficiency
of the cell problem solution can be significantly improved using
\gls{fft}-based methods~\cite{Kochmann2018, Roters2013, Prakash2009}.

As fast and reliable \gls{rve} solvers, \gls{fft}-based methods have gained
much attention in the last 20 years~\cite{Vondrejc2014,Zeman2017,
  Schneider2021}.~\gls{fft}-based methods exploit the simple structure
of regular grids and allow one to use lightweight iterative solvers such as
Newton-\gls{cg} for solving computational homogenization
problems~\cite{Zeman2010, Mishra2016}. By contrast, in the conventional
\gls{fem} framework, due to poor scaling of linear solution complexity (caused
by deterioration of spectral characteristics of the corresponding linearized
system with growing the size of the problem), iterative solvers are not a good
choice unless being equipped with preconditioners and direct solvers are the
used most of the time~\cite{Moshfegh2020}. \gls{fft}-based methods are,
therefore,  roughly 200 times faster compared to the conventional~\gls{fem}
scheme~\cite{Eisenlohr2013} solving problems with roughly $10^4$~\gls{dof}s and
higher.
\gls{fft}-based methods are also less memory consuming (due to their inherent
matrix-free formulation~\cite{Schneider2021}) compared to~\gls{fem}. In
addition, solution domain discretization is trivial since the~\gls{rve} is only
implicitly meshed.

Increasing the resolution and improving solution accuracy is straightforward in
spectral methods.
As a result, large-scale micro-structures simulations are more efficient using
spectral methods compared to conventional~\gls{fem}.
The computational complexity of~\gls{fft}-based solution
in the \gls{fg} scheme~\cite{deGeus2017} is dominated by the~\gls{fft} algorithm
($\mathcal{O}\left(n \log{n}\right)$, where $n$ is the number of the
discretization points in the \gls{rve} solution domain).
The availability of highly optimized \gls{fft} implementations
(\textit{FFTW}~\cite{Frigo2005a} and \textit{PFFT}~\cite{Pippig2013}) makes
efficient implementation of spectral methods simple.

Spectral methods, in addition to their original scheme (iterative solution
of an integral equation of the Lippmann-Schwinger type) as introduced by
Moulinec et al. \cite{Moulinec1994, Pivovarov2018} can be derived by applying
the Galerkin method using trigonometric polynomials as shape functions
(\gls{fg} method)~\cite{Vondrejc2014, Zeman2017}.
The global support and the oscillatory nature of trigonometric functions
results in Gibbs ringing phenomenon~\cite{Ma2021a} which makes them unsuitable for
problems containing localized phenomena such as continuum damage and non-liner
plasticity.
Gibbs ringing artifact occurs near abrupt transitions as spurious fluctuations
in the solution of the problem (strain and stress field in case of mechanical
computational homogenization) due to discrete Fourier transformation truncation.
Gibbs-ringing is well-known and well-documented, for instance
in~\cite{Gottlieb1997, Hewitt1979, Gelb2007, Gottlieb1996, Ma2021a}.
Sharp discontinuity with high phase contrast
exacerbates Gibbs-ringing, thus, Gibbs ringing is more pronounced in problems
containing phases with high contrast, for instance
\gls{rve}s containing highly contrasted composites~\cite{Michel2001} or
void in the micro-structure such as foams~\cite{Marvi-Mashhadi2020}.

Several approaches have been used to address Gibbs-ringing
in spectral methods, e.g. ~\cite{Willot2015, Schneider2016, Khorrami2020,Ma2021}.
These methods can be categorized into two main categories
\rom{1}. mitigation: these approaches try to reduce the
Gibbs-ringing artifact fluctuations by smoothing the geometry of the phase
interfaces~\cite{Kabel2015} or by using higher order discrete
derivatives~\cite{Schneider2016}, or \rom{2}. elimination: removing the
Gibbs-ringing by using discrete derivatives obtained from a regular periodic
\gls{fem} discretization~\cite{Leute2022, Leuschner2018}. The projection based method
in~\cite{Leute2022} is capable of eliminating the Gibbs ringing while
maintaining all the advantages of the~\gls{fft}-based methods mentioned
earlier. A preconditioned
displacement based~\gls{fem} scheme mathematically equivalent to the
projection-based scheme of~\citet{Leute2022} has been developed recently
by~\citet{Ladecky2022}. In their FFT-accelerated \gls{fe} scheme, a
preconditioner derived based on the Green's function of a uniform
periodic reference medium is used to make the distribution of system matrix
eigenvalues favorable for iterative solvers~\cite{Ladecky2021a, Pultarova2021}.
Such discretized Green's functions are mathematical convolutions and therefore
their inversion and application is cost-effective in the Fourier space.

Use of Newton-\gls{cg} due to its quadratic convergence properties is
common in \gls{fft}-based solution schemes~\cite{Gelebart2013, Kabel2014,
  Vinogradov2008, Zeman2017}.
However, it is notable that Newton-\gls{cg} is unable to
handle problems with non-convex energy functional (\gls{snpsd} system matrix)
such as homogenization of \gls{rve}s containing meta-materials~\cite{Ai2017} or
continuum damage~\cite{Bazant1976, Marvi-Mashhadi2020}.
Conventional non-linear \gls{fem} solvers are also susceptible to instabilities
in modeling non-convex problems~\cite{Pijaudier-Cabot1987}.

A possible solution, called \gls{sla} ~\cite{Rots2001, Rots2004, Rots2008,
  Dejong2008, Rots1988}, for circumventing the
solution of non-convex problems is breaking the non-linear non-convex problem
to a sequence of linear convex problems. In~\gls{sla}, in each solution step,
only one integration point is allowed to soften by certain reduction of
its stiffness due to damage~\cite{Pari2022}. \gls{sla} is an event-driven algorithm and
therefore it does not scale with problem size since by increasing problem size
damage sites (possible event sites) increases rapidly; therefore, \gls{sla}
becomes inefficient. This is a major drawback for solving a cell problem in a
multi-scale model~\cite{Ramos2017}.

Non-convex iterative solvers, for instance nonlinear~\gls{cg},
quasi-Newton solvers such as \gls{lbfgs}~\cite{Nocedal2006,
  Curtis2015}, and also trust region solvers~\cite{Yuan2000,Nocedal2006}
are typically developed for optimization of problems with known objective
functions. However, in mechanical engineering problems, the objective function
(strain energy functional) is often not known or difficult to calculate for
complicated constitutive laws. In context of~\gls{fft}-accelerated solvers
~\cite{Wicht2019} has used quasi-Newton solvers for convex problems using an
inexact line search method, but the application to non-convex problems still
remains not investigated.

In this paper, we introduce a~\gls{faief}
that enables application of non-convex solvers for problems whose
objective function explicit form is not available.
The introduced~\gls{faief} is a first order
approximation and is only valid for small load steps.
However, this is typically not an issue because many homogenization problems, such as
nonlinear computational homogenization problems, already meet this condition
as for solving them the load increments are chosen to be small.

We opted for trust-region Newton-\gls{cg} solver
to illustrate the applicability of the introduced~\gls{faief} for solving
non-convex problems whose energy functional evaluation is not easy. In addition,
we demonstrate how the projection-based scheme with~\gls{fe} discretization
~\cite{Leute2022}, in conjunction with the introduced trust region
solver~\cite{Nocedal2006} (equipped with the introduced~\gls{faief}),
can successfully and efficiently solve homogenization damage problems
with a ~\gls{snpsd} Hessian matrix.

The projection-based scheme with~\gls{fe} discretization as well as the
employed non-convex solvers are explained in Section~\ref{chap:methods}.
In Section~\ref{chap:Results}, we show that the
modified trust region Newton-\gls{cg}, developed and introduced in this paper,
yields identical results compared to the standard trust region Newton-\gls{cg}
solver in a simplistic $1$D example. The convergence properties of
the modified trust-region solver are resolution independent, as we show by
presenting the modified trust region Newton-\gls{cg} solver performance for
variation of \gls{rve} sizes. In addition, we show that the  modified trust
region solver can solve a real-world damage homogenization problem that suffers
both from \gls{snpsd} stiffness and sharp phase interface with unbounded
contrast.
We have implemented all methods in the open source code
$\mu$Spectre~\cite{muspectre} and all numerical examples in
Section~\ref{chap:Results} can be reproduced by executing the
corresponding scripts provided in the supplementary material.

%%%%%%%%%%%%%%%%%%%%%%%%%%%%%%%%%%%%%%%%%%%%%%%%%%%%%%%%%%%%%%%%%%%%%%%%%%% 
%                                                                         % 
% Methods                                     %
%                                                                         % 
%%%%%%%%%%%%%%%%%%%%%%%%%%%%%%%%%%%%%%%%%%%%%%%%%%%%%%%%%%%%%%%%%%%%%%%%%%% 
\section{Methods}
\label{chap:methods}
\begin{figure}
  \centering
  \includegraphics[height=6.0cm, width=6.0cm]{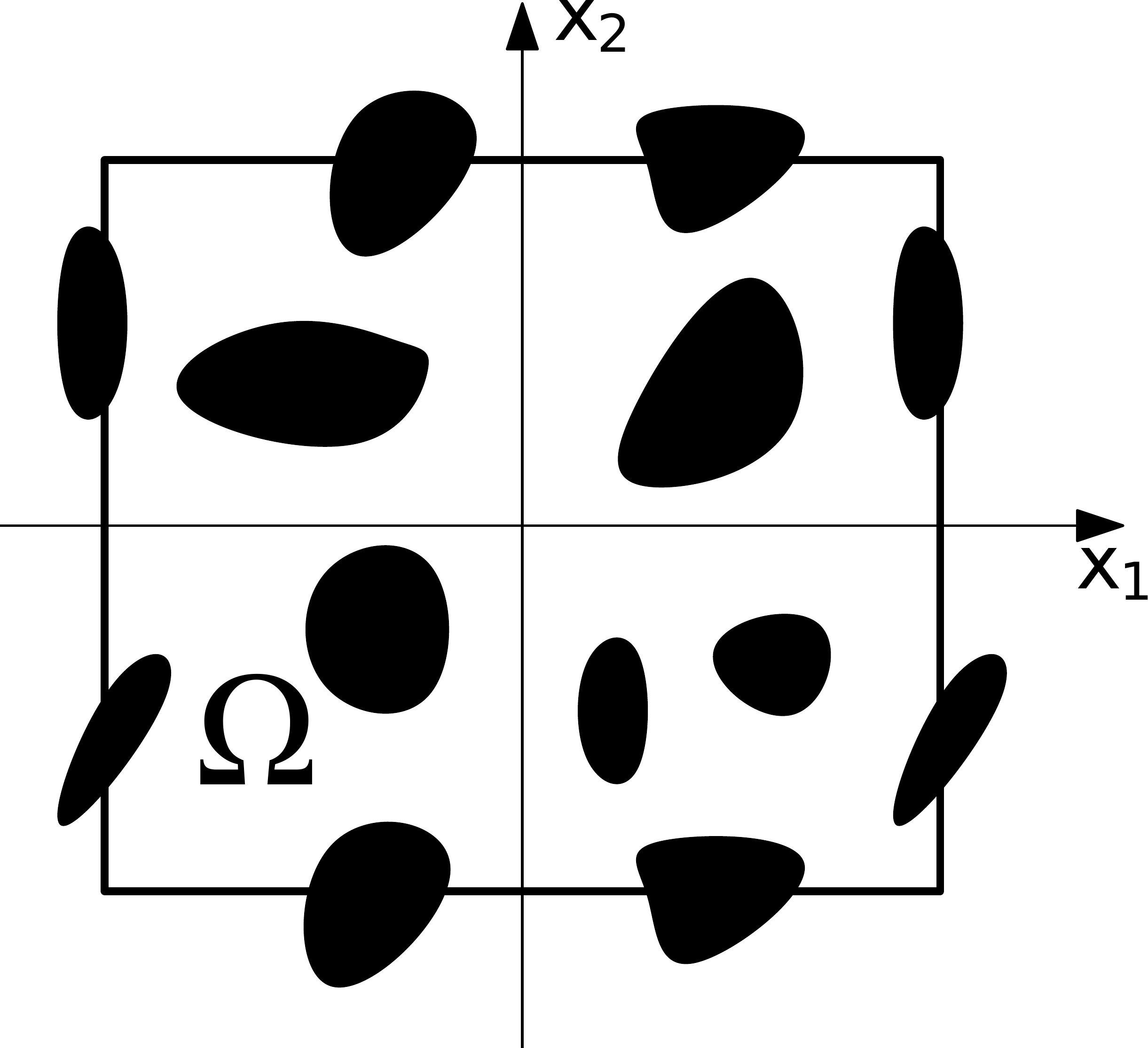}
  \caption{Rectangular two-dimensional cell with outlined periodic
    micro-structure obtained from Ref.~\cite{Ladecky2022}.}
  \label{fig:illustration}
\end{figure}

In the following, we consider a rectangular periodic \gls{rve};
see~\autoref{fig:illustration} for illustration of a typical micro-structure.
Small strain micro-mechanical formulation is adopted for the derivations of the
equilibrium of a
micro-structure undergoing a displacement field of $\vek{\chi}: \pucc\rightarrow\Omega$ which maps
the grid points from undeformed positions $\pucc$ to their deformed
configurations $\Omega$. The material response
corresponding to position $\vek{x}$, given the local strain $\strain(\vek{x})$
is determined by constitutive law of the material at that point as $\stress
\left( \vek{x}, \strain(\vek{x}) \right)$.
The total strain $\strain$ can be decomposed to average strain tensor $\vek{E}$
and the periodic fluctuating strain field $\fluc{\vek{\strain}}(\vek{x})$:
\begin{align}
  \label{eq:strain_decomp}
  \strain(\vek{x}) = \vek{E} + \fluc{\vek{\strain}}(\vek{x})\ \forall \vek{x} \in \pucc,&& \text{s.t.} &&
                                                                            \int_{\pucc}\fluc{\strain}(\vek{x})
                                                                            d\vek{x} = \vek{0}
\end{align}

The governing mechanical equilibrium equation in this
domain reads as:
\begin{subequations}
  \label{eq:equilibrium}
  \begin{align}
    \label{eq:equilibrium1}
    -\nabla_0 \cdot \vek{\stress} \left(\vek{E} + \fluc{\vek{\strain}}(\vek{x}),
    \vek{x}\right)&=0,\ \forall \vek{x} \in \pucc.
  \end{align}
  \begin{align}
    \label{eq:compatibility}
    \fluc{\vek{\strain}} &= \nabla_{0, s}\fluc{\vek{u}} \in \mathcal{E} && \fluc{\vek{u}}:\pucc\
                                                                    \text{periodic
                                                                    displacement}
  \end{align}
\end{subequations}
where $\nabla_0$ is the \textit{nabla} operator ($\nabla_0$ is gradient operator
and $\nabla_0\cdot$ is the divergence operator) in the reference coordinates and
$\nabla_{0,s}$ stands for the symmetrized \textit{nabla} operator
$\left(\nabla_{0,s} = \frac{1}{2}(\nabla_0+ \nabla_0^T)\right)$ in the reference
coordinates. The compatibility equation \eqref{eq:compatibility}
dictates the strain field to be a gradient field.

In the projection-based spectral scheme, in order to solve the equilibrium
equation, the local equilibrium equation in its weak form after applying Gauss
divergence theorem casts to:

\begin{equation}
  \label{eq:weak_form}
  \int_{\Omega_{0}} \compat{\test}(\vek{x}):\stress \left(\vek{E} +
    \fluc{\vek{\strain}}(\vek{x}), \vek{x}\right) d\vek{x} = 0 \qquad \forall \compat{\test} \in
  \mathcal{E}.
\end{equation}
Note that the boundary term cancels out due to periodic boundary condition. The
strain is taken to be from the compatible gradient space $\mathcal{E}$.

The test function in the projection-based spectral methods ($\compat{\test}$ in
\eqref{eq:weak_form}) is in strain space
while in the \gls{fem} formulation, the test function is a displacement field.
Therefore, the test function of projection-based spectral methods is not an
arbitrary field (unlike \gls{fem}). Proceeding with the Galerkin discretization
necessitates having a fully arbitrary test function. Therefore, in order to
apply the Galerkin discretization on \eqref{eq:weak_form},
it is necessary to impose compatibility condition on the test variable
$\test(\vek{x})$.~\citet{Zeman2017} have introduced a compatibility
projection operator $\mathbb{G}$ based on the Fourier discretization to impose
compatibility. The operator $\mathbb{G}$ maps any second-order tensor to
its compatible (periodic gradient) contribution. Applying the projection
operator on the test variable  $\test(\vek{x})$ makes it possible to
continue with the Galerkin discretization and solve directly for strain field.

In the following of this section, after elaborating the
projection operator $\mathbb{G}$ and the resulting discretized equilibrium
equation, the remedy of the ringing artifact developed by Leute
et. al~\cite{Leute2022} is reviewed. Next, its extension to
trust region solver capable of handling generic homogenization problems with
non-convex energy functions such as damage mechanics problems is explained.

%%%%%%%%%%%%%%%%%%%%%%%%%%%%% projection based solvers     %%%%%%%%%%%%%%%%%%
\subsection{Projection based solver}
\label{sec:projection_based_solver}
% Spectral solvers (Richard Paper)
The key element of projection-based spectral solvers is the projection operator
$\mathbb{G}$ which enforces compatibility on an arbitrary field
$\test(\vek{x})$ as:

\begin{align}
  \label{eq:projection}
  \compat{\test}(\vek{x}) &= \left[\mathbb{G} \star  \test\right](\vek{x}) \nonumber \\
                          &=\int_{\pucc}\mathbb{G}(\vek{x}-\vek{y}):\test(\vek{y})\
                            d\vek{y} \ \  \forall \vek{x} \in \pucc.
\end{align}
Applying projection yields compatible contribution $\compat{\test}(\vek{x})$ of
the original field $\test(\vek{x})$.
In~\eqref{eq:projection}, $\star$ is the convolution operator. The convolution
format of~\eqref{eq:projection} makes its application in
Fourier space convenient, since convolution in real space is equivalent to
contraction in Fourier space. Accordingly, \eqref{eq:projection} can be
rewritten as:
\begin{equation}
  \label{eq:conv_to_fourier}
  \compat{\test}(\vek{x}) =
  \iFourierTrans\{\Fouri{\mathbb{G}}(\vek{k}):\Fouri{\zeta}(\vek{k})\}
\end{equation}
where $\Fouri{\mathbb{G}}(\vek{k})$ is the compatibility operator in Fourier basis,
$\Fouri{\zeta}(\vek{k})=\FourierTrans \{ \zeta( \vek{x} ) \}$, and
$\vek{k}$ is the discrete frequency vector in the Fourier domain.
Considering that $\Fouri{\mathbb{G}}$ is a fourth order block diagonal tensor,
the equilibrium equation~\eqref{eq:weak_form} combined with the compatibility
equation~\eqref{eq:conv_to_fourier} yields~\cite{Zeman2017}:
\begin{equation}
  \label{eq:project_stress}
  \Fouri{\mathbb{G}}(\vek{k}):\Fouri{\stress}(\vek{k})=\vek{0}.
\end{equation}

Using Newton's method to solve the nonlinear system of \eqref{eq:project_stress}
iteratively, the $(i+1)^{th}$ update of the strain field
$\vek{\strain}_{i+1}$ in the iterative scheme can be calculated from the
previous approximation of the strain field $\vek{\strain}_{i}$ incremented by
a finite strain increment $\delta \vek{\strain}_{i+1}$,
\begin{equation}
  \label{eq:strain_linearization}
  \vek{\strain}_{i+1} = \vek{\strain}_{i} + \delta \vek{\strain}_{i+1}.
\end{equation}
Starting from an initial strain approximation $\vek{\strain}_{0}$ the strain
increment at each step is given by solution of the linear system:
\begin{equation}
  \label{eq:project_stress_linear}
  \Fouri{\mathbb{G}} : \tangent_{i} : \delta \boldsymbol{\strain}_{i+1} =-\
  \Fouri{\mathbb{G}} :
  \boldsymbol{\stress}_{i}
\end{equation}
at the $(i+1)^{th}$ nonlinear solution step, where $\tangent_{i}$ is the
constitutive tangent matrix of the system evaluated at discretization points.

\citet{Leute2022} derived a general expression for the operator
$\Fouri{\mathbb{G}}$ as an explicit function of the second rank tensor
$\Fouri{\vek{g}}$ with the form of:
\begin{align}
  \label{eq:proj_operator}
  \Fouri{\vek{g}}(\vek{k}) = \frac{\Fouri{\vek{\Doper}} (\vek{k})\otimes
  \HTranspose{\Fouri{\vek{\Doper}}}(\vek{k})}{\Fouri{\vek{\Doper}}(\vek{k})\cdot
  \HTranspose{\Fouri{\vek{\Doper}}}(\vek{k})},
\end{align}
where $\Fouri{\vek{\Doper}}$ is the derivative operator in Fourier space, and $^*$
denotes the Hermitian transpose. The form of the projection operator
$\Fouri{\mathbb{G}}$ as a function of $\Fouri{\vek{g}}$ is different in small
strain and finite strain formulation (further details can be found
in~\cite{Zeman2017, Leute2022}).

In the original projection based method developed by~\citet{Zeman2017},
$\Fouri{\vek{\Doper}}$ was expressed based on the Fourier basis as
$\Fouri{\Doper}(\vek{k})=i \vek{k}$ which yields a second rank tensor
$\Fouri{\vek{g}}$ of the form:
\begin{align}
  \label{eq:proj_Fourier_2}
  \Fouri{\vek{g}}_{\alpha\beta}(\vek{k})
  =  \begin{cases}
    \t{0} & \mathrm{if}\,\ \vek{k}=\vek{0} , \\
    \frac{\vek{k}_\alpha \vek{k}_\beta}{\vek{k}^2} & \forall \  \vek{k} \neq \vek{0},
  \end{cases}
\end{align}
in index notation, where $\vek{k}$s are normalized discrete Fourier wave-vectors.
\citet{Leute2022} showed that, based on the general form of the projection
operator ($\Fouri{\mathbb{G}} =
\Fouri{\mathbb{G}}\left(\Fouri{\vek{g}}\right)$ explicitly formulated
in~\cite{deGeus2017} in both finite-strain and small strain formulation), it is
possible to derive projection operators using gradient operator obtained from
different discretization schemes of choice. For instance, \citet{Leute2022}
worked out a projection operator based on a linear \gls{fe} discretization and
showed that using the basis functions of a \gls{fe} discretization
results in elimination of Gibbs ringing artifacts.
This makes use of \gls{fe} discretization suitable for
problems with highly localized phenomena such as damage mechanics. In the
following of this paper, the derivations are carried out for a spectral method
with \gls{fe} basis set.

Choosing different $\Fouri{\Doper}$ in a projection based spectral method is
equivalent to choosing different element types and shape functions in the
conventional \gls{fem} formulation. The operator $\Fouri{\Doper}$, for
a~\gls{fe} discretization, is calculated using the derivative of
the corresponding shape functions.
Similar to a \gls{fe} scheme, the strain, stress, and the constitutive tangent are
evaluated at the quadrature point of the \gls{fe} discretization.

Optimal spectral characteristics of the system matrix in
\eqref{eq:project_stress_linear}~\cite{Pultarova2021}
makes linear iterative solvers and specifically \gls{cg} solver ideal for
solving it. However, solving~\eqref{eq:project_stress_linear} with \gls{cg}
solver needs Hessian matrix $\tangent_{i}$ to be \gls{spsd} which is not the
case in several mechanical homogenization problems such as
meta-materials~\cite{Ai2017} or continuum damage~\cite{Bazant1976,
  Marvi-Mashhadi2020}. In the~\autoref{sec:non_convex_solver}, a trust region
Newton-\gls{cg} solver is adopted to expand the use of the projection-based
spectral method to non-convex problems.

%%%%%%%%%%%%%%%%%%%%%% Non-convex solvers     %%%%%%%%%%%%%%%%%%%%%%%%
\subsection{Non-convex optimization}
\label{sec:non_convex_solver}
The projection-based spectral scheme explained in
Subsection~\ref{sec:projection_based_solver} due to the optimal spectral
characteristics of its linearized system~\cite{Pultarova2021}
enables us to benefit from the
computational advantages of iterative solvers such as Newton-CG. This results in
great scaling for solving \gls{rve}s with large number of discretization
points. As set forth above, \gls{cg} can only solve linear system of equations
with \gls{spsd} matrix.

Therefore, in order to be able to benefit from the computational
speed-up offered by the projection based spectral scheme for problems with
\gls{snpsd} system matrices, we need to employ other iterative solvers.
Possible candidates capable of
handling non-convex problems are nonlinear~\gls{cg}, quasi-Newton solvers, and
Trust region Newton solver. In this paper, a trust region Newton solver as
a robust and memory efficient solver capable of handling non-convexity in an
iterative fashion is adopted to be used in homogenization problems.

The potential of trust region solvers (as well as quasi-Newton solvers)
has not been exploited in computational homogenization, since in a considerable
part of the literature, conventional \gls{fem} employing direct solvers used to
be the de-facto for computational homogenization. As a result,
the main challenge of using trust region and quasi Newton solvers
in computational homogenization, namely missing an explicit expression of the
objective function in the equivalent energy minimization counterpart of the
equilibrium solution, has not been addressed to the authors' best
knowledge. This problem is addressed in the following of this section after a
review over the standard trust region Newton solver.

\subsubsection{Trust region Newton solver}
In contrast to conventional line search algorithms, in a trust region solver,
the approximate model (sub-problem) "trusted" within a bounded region (trust
region) near the current iterate is minimized iteratively until a minimizer of
the original function is reached~\cite{Nocedal2006, Byrd2000,
  Conn2000}.
In the trust region solution strategy, optimization is carried out by
minimizing a model function (typically quadratic) trusted up to a certain
radius of the current iterate as a proxy problem (sub-problem). The approximate
model is derived from the local information gathered from the objective
function at the current iterate.

There is no a-priori knowledge of the radius of the region in which the model
can adequately approximate the objective function. In addition, the model's
accuracy declines by moving away from the current iteration. Therefore, it
is crucial to determine the "trust region" of the model function and to
regulate it consistently during the solution process.

To this end, based on the model's match with the original objective
function, the trust region is adjusted in each iteration. As a
general rule, the trust region can be expanded if the approximate model fits
the original problem well. In contrast, the trust region shrinks if the
approximate model fails to estimate the original function
adequately~\cite{Hsia2017}. It is therefore necessary
to access the original objective function of the equivalent optimization
problem if one needs to use a trust region solver in its standard formulation.

Let us consider the total energy function, equivalent to the original
objective function of a homogenization problem, of the \gls{rve} as:
\begin{equation}
  \label{eq:energy_func}
  W=
  \sum_{Q}{w(\vek{\strain}^Q,\ \vek{g(\vek{x}^Q)})}
\end{equation}
where $w$ denotes energy at $Q$s, which are the discretization
quadrature points, and $\vek{g(\vek{x})}$ represents the internal variables of
the material.

The solution of the equilibrium equation \eqref{eq:equilibrium1} is,
in particular, corresponding to the critical point of the total energy
function \eqref{eq:energy_func}.
In the projection-based scheme the equilibrium problem is solved in
strain state space. According to the fact that the energy conjugate of strain
is stress, the gradient of the energy with respect to strain is actually the
stress tensor. In addition, the Hessian of the objective function
corresponds to the tangent stiffness of the material at the discretization points,
\begin{align}
  \label{eq:graient_and_hessian}
  % \nablastrain W &= \frac{\partial W}{\partial \vek{\strain}} = \vek{\stress} \\
  % \nablastrain^2 W &= \frac{\partial^2 W}{\partial \vek{\strain}^2} = \tangent
  \nablastrain W   &= \vek{\stress}, \\
  \nablastrain^2 W &= \tangent.
\end{align}
% we might need to mention hyper-elastic assumption
One possible sub-problem model (and probably the most common form)
of the trust region solver is a quadratic energy function approximation in form
of:
\begin{equation}
  \label{eq:modeling_trust}
  m_i(\vek{p}_i) = W(\vek{\strain}_i) + \nablastrain W^T\vek{p}_i +
  \frac{1}{2}\vek{p}_i^T\tangent_i\vek{p}_i, \qquad  s.t.\ ||\vek{p}_i||<R_i,
\end{equation}
where in the Newton trust region solver,
$\tangent_i$ is taken as the Hessian matrix of
the energy evaluated $\left(\tangent \right)$ at quadrature points at
$i^{\text{th}}$ load step and $R_i$ is the radius of the trust region.
$\vek{p}_i$ is a solution step in the strain space.
Here, the energy functional is taken as a direct function of strains at
$i^{\text{th}}$ solution step as in the projection based formulation the
equations are solved in strain space. Other trust region solvers are also
possible using different choices for the matrix $\tangent_i$.

The agreement of the actual objective function
($W(\vek{\strain}_i + \vek{p}_i)$) and the model ($ m_i(\vek{p}_i)$) at the new
iterate is evaluated by a scalar variable at the $i^{\text{th}}$ iterative
step, defined as:
\begin{equation}
  \label{eq:rho}
  \rho_i = \frac{W(\vek{\strain}_i) - W(\vek{\strain}_i + \vek{p}_i)}{m_i(\vek{0})
    - m_i(\vek{p}_i)}.
\end{equation}

In the trust region algorithm, the value of  $\rho_i$, as set forth
in~\cite{Nocedal2006}, determines how the trust region size will be
updated as well as whether or not the proposed step will be accepted.

It is relatively simple to calculate the denominator of the right hand side of
the~\eqref{eq:rho} ($\Delta m_i =m_i(\vek{0})- m_i(\vek{p}_i)$) according to the
definition of $m_i(\vek{p}_i)$ given in~\eqref{eq:modeling_trust}.
For calculation of the nominator, explicit expression of the origin
objective function (equivalent to stored energy in mechanical engineering
problems) is necessary. In mechanical homogenization problems, however, the
objective function is often not calculable (at least easily), since the actual
energy density function of most materials is very complex (and even impossible
to compute). A first order incremental energy function that allows us to use
trust region solvers to solve generic non-convex mechanical homogenization
problems is presented here. In order to derive the incremental energy
functional, first, Taylor series of the actual energy function $W$ is expanded
at both $\vek{\strain}_i$ and $\vek{\strain}_i+\vek{p}_i$ points as:
\begin{subequations}
  \begin{equation}
    \label{eq:Taylor1}
    W(\vek{\strain}_i + \vek{p}_i) = W(\vek{\strain}_i) + (\nablastrain
    W|_{\vek{\strain}_i})^T \vek{p}_i + \frac{1}{2} \vek{p}_i^T (\nablastrain^2
    W|_{\vek{\strain}_i}) \vek{p}_i + \higherorders,
  \end{equation}
  \begin{align}
    \begin{split}
      \label{eq:Taylor2}
      W(\vek{\strain}_i ) &= W((\vek{\strain}_i + \vek{p}_i) - \vek{p}_i)\\
      &= W(\vek{\strain}_i + \vek{p}_i) - (\nablastrain W|_{\vek{\strain}_i+\vek{p}_i})^T
      \vek{p}_i + \frac{1}{2} \vek{p}_i^T (\nablastrain^2 W|_{\vek{\strain}_i+\vek{p}_i})
      \vek{p}_i + \higherorders.
    \end{split}
  \end{align}
\end{subequations}

Subtracting \eqref{eq:Taylor1} from \eqref{eq:Taylor2} and dropping higher
order terms yields:
\begin{align}
  \label{eq:energy_estim}
  W(\vek{\strain_i} + \vek{p_i}) - W(\vek{\strain_i} ) &\approx\\ \nonumber
                                                       &\frac{1}{2}\left((\nablastrain
                                                         W|_{\vek{\strain_i}})^T_i
                                                         \vek{p}_i
                                                         + (\nablastrain
                                                         W|_{\vek{\strain_i}+\vek{p_i}})^T_i
                                                         \vek{p}_i\right)+\\ \nonumber
                                                       &\frac{1}{4}\left(\vek{p}_i^T (\nablastrain^2 W|_{\vek{\strain_i}})\vek{p}_i -
                                                         \vek{p}_i^T (\nablastrain^2 W|_{\vek{\strain_i}+\vek{p_i}}) \vek{p}_i\right)
\end{align}
Truncating \eqref{eq:energy_estim} up to first order gives:
\begin{align}
  \label{eq:incr_energy_func}
  W(\vek{\strain_i} + \vek{p_i}) - W(\vek{\strain_i} ) \approx \overline{\Delta W} =
  \frac{\stress(\vek{\strain}_i+\vek{p}_i)+\stress(\vek{\strain}_i)}{2}:\vek{p}_i.
\end{align}
The terms in the right hand side of
~\eqref{eq:incr_energy_func} consist of the stress tensors at the previous
and current trial steps which are already evaluated at all
of the quadrature points.
To take $\overline{\Delta W}$ as a valid estimation of $\Delta W$, it is
necessary to keep the load increments small.
Furthermore, it is vital that the variation of resulting displacement field
remains bounded. For instance, problems such as buckling under prescribed
growing force does not satisfy boundedness of displacement field around the
critical load; and therefore does not converge using the modified trust-region
solver presented here.
On the other hand, problems such as mechanics damage modeling are solvable
using the modified solver presented here given that the applied load increments
are controlled to be small.

The actual system energy reduction in the nominator of~\eqref{eq:rho} can be
replaced by the
first order incremental energy approximation $\overline{\Delta W}$ calculated by
energy~\eqref{eq:energy_estim} which gives an estimation of $\rho$; denoted by
$\overline{\rho}$ according to:
\begin{equation}
  \label{eq:rho_bar}
  \bar{\rho}_i = \frac{\overline{\Delta W}}{\Delta m_i}.
\end{equation}
It is notable that the evaluation of $\overline{\rho}_i$ needs the stress of
the previous solution step to be stored.

Introducing the first order approximation of the energy functional
enables us to use the robust trust region algorithm in cases that the explicit
expression of $W$ function is missing. The pseudo-algorithm
of the projection-based trust-region solver is presented
in~\autoref{algo:Trust region Newton-CG} which shows how the scalar value
$\rho$ (or its estimated counterpart $\overline{\rho}_i$) is used to
make decisions of accepting or rejecting trial step ($\delta \t{F}$) as well as
shrinking or expanding the trust region. As shown in~\autoref{algo:Trust region
  Newton-CG}, the memory overhead of using Trust region Newton-\gls{cg} is
storing the flux (stress) field at the previous solution step which does not impact the
overall required for the solution compared to Newton-\gls{cg}.

The predicted reduction of the model ($m_i$) will always be non-negative since
the step $\vek{p}_i$ is calculated by minimizing the model $m_i$ over the
region that includes $\vek{p=0}$. Therefore, if $\overline{\rho}_i$ is negative,
the value of the objective function at the new iterate
$\left(W(\vek{\varepsilon}_i + \vek{p}_i)\right)$ is greater than the current objective value $\left(W
  (\vek{\varepsilon}_i)\right)$, thus the step must be rejected. Alternatively, if
$\overline{\rho}_i$ is close to 1, it is safe to expand the trust
region for the next step since the model $m_i$ and the original objective
function $W$ are in good agreement over the
solution step $\vek{p}_i$. In the case when $\overline{\rho_i}$ is positive but smaller than
one, the trust region is not altered, however in the case when $\overline{\rho}_i$
is close to zero or negative, the trust region is shrinked by reducing $R_i$ at
the next iteration. How these decisions are made is depicted in details in
~\autoref{algo:Trust region Newton-CG}.

The solution of the subproblem (minimizing $m$ within the ball of radius $R$) is
easy to find when $\tangent$ is positive definite and the minimizer is
located within the trust region (equivalent to finding an unconstrained
minimizer of the quadratic function $m$).
There is no such simple solution to the subproblem in other cases.
For these cases
the minimzer resides on the boundary of the trust region. The constrained
linear solver used here as the sub-problem solver is based on the solver
introduced by Steihaug~\cite{Nocedal2006, Steihaug1983} which is used to
generate the trial solution step $\vek{p_i}$.

In addition, in order to make the
linear solver robust to numerical rounding error problems a reset algorithm
based on the work of~\citet{Powell1977}, and~\citet{Dai2004} was added to the
linear solver algorithm. The used reset mechanism
replaces the conjugate gradient step with a restart step (for instance, the
gradient descent step) whenever two successive solution steps inside the linear
solver ($\vek{r}_j$ and $\vek{r}_{j-1}$) are not
sufficiently orthogonal to each other. The measure expressing the orthogonality
of the solution steps can be calculated by:
\begin{equation}
  \label{eq:1}
  r_{\sphericalangle} = \frac{|\vek{r}_j \cdot \vek{r}_{j-1}|}{||\vek{r}_j^2 ||}.
\end{equation}
Comparing the measure $r_{\sphericalangle}$ by a constant value in the range of $(0.1, 0.9)$ has
been proposed as the decision criteria for restarting \gls{cg}. We chose
$0.2$ as suggested in~\cite{Powell1977}, hence the restart procedure is invoked
if the measure $r_{\sphericalangle}$ is greater than $0.2$.

Several quasi-Newton solvers also depend on the explicit
expression of the objective function (e.g. \gls{lbfgs}).
Using the approximated strain energy reduction $\overline{\Delta W}$
instead of $\Delta W$ makes use of these quasi-Newton solvers possible as well.
As derived here, the introduced incremental approximation of the objective
function, specifically,
enables us to use trust region Newton-CG solver following the
algorithm given in~\autoref{algo:Trust region Newton-CG} (as the Newton
nonlinear solver) and the algorithm given in~\autoref{algo:CG-Steihaug solver}
as the subproblem solver.  The introduced approximate energy
functional can be generalized to solve any other problem in which the explicit
objective function is not available or not easy to calculate while the gradient
and Hessian of the objective function are available.

%%%%%%%%%%%%%%%%%%%%%%%%%%%%%%%%%%%%%%%%%%%%%%%%%%%%%%%%%%%%%%%%%%%%%%%%%%%% 
%                                                                          % 
% Results                                                                  %
%                                                                          % 
%%%%%%%%%%%%%%%%%%%%%%%%%%%%%%%%%%%%%%%%%%%%%%%%%%%%%%%%%%%%%%%%%%%%%%%%%%%% 
\section{Results and discussion}
\label{chap:Results}
In the following, three examples  are presented to demonstrate the solver
developed above are presented. We first compare the performance of the modified
trust region solver with that of the Newton-CG solver and the standard trust
region solver on a very simple example. This example deals with a $1$D damage
spring system where the potential energy functional of the system is known. The
trust-region solver can therefore be used on this system. Second, we test the
correctness of the solution of the solver considering a convex system. The
Eshelby inhomogeneity, whose analytical solution is available, is selected, so
we can compare the generated solution of the modified trust region solver with
the analytical solution of the problem. Finally, the capability of the solver
for solving non-convex problems is demonstrated in a real world damage
mechanics problem. As illustrated in this example the modified trust region
solver can handle damage mechanics homogenization problems with rather complex
constitutive laws.

\subsection{Minimal 1D non-convex example}
\label{sec:1d_example}
\begin{figure}
  \captionsetup{}
  \centering
  \begin{subfigure}{.5\textwidth}
    \centering
    \vspace{0.75cm}
    \includegraphics[height=2.0cm, width =
    5.0cm]{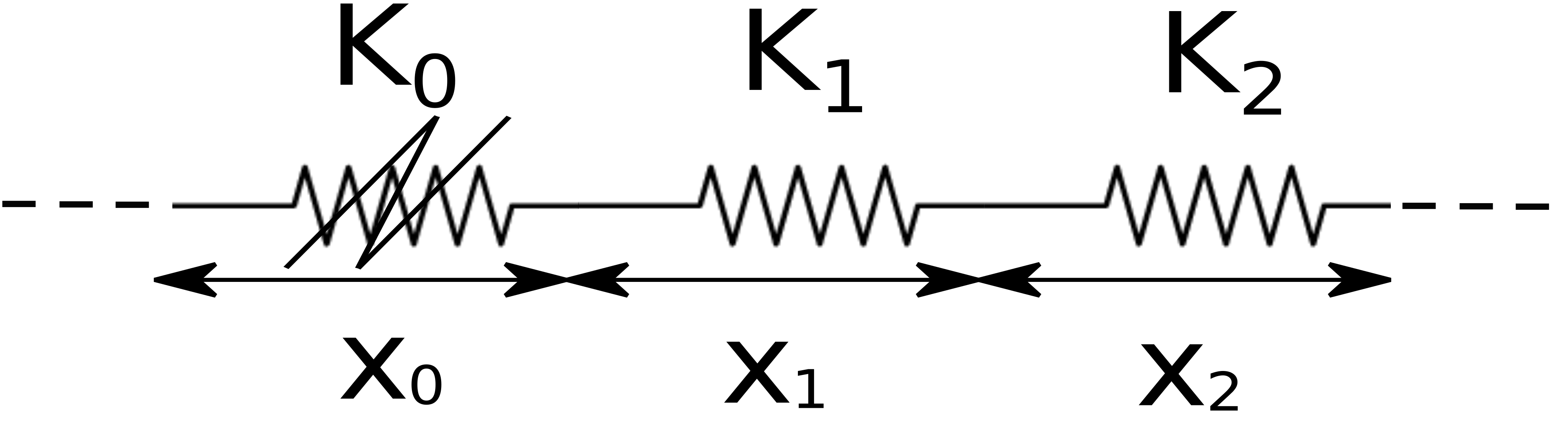}
    \vspace{0.75cm}
    \caption{ }
    \label{fig:spring_schematic}
  \end{subfigure}%
  \begin{subfigure}{.5\textwidth}

    \centering
    \includegraphics[height=4.5cm, width =
    5.0cm]{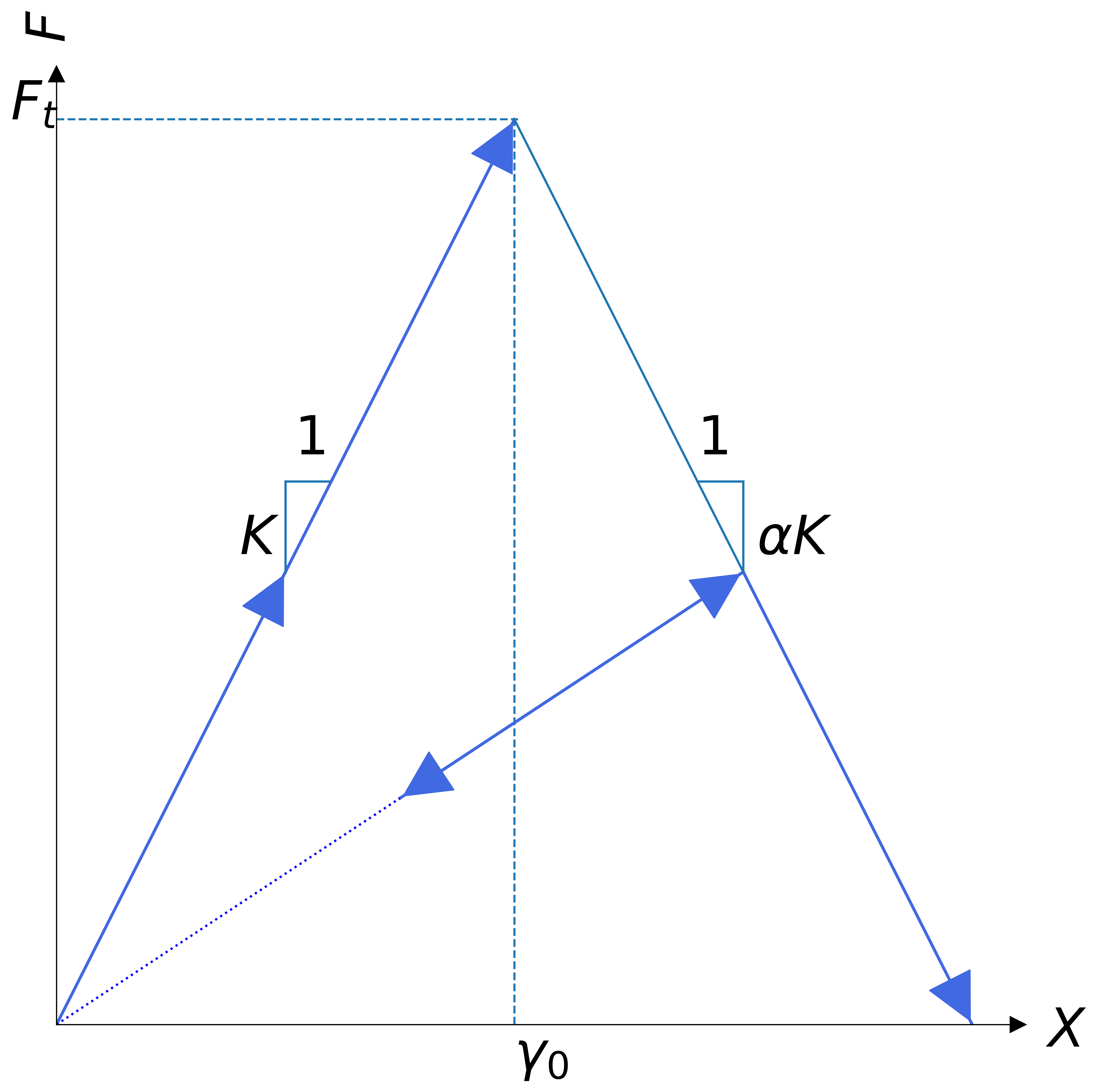}
    \caption{ }
    \label{fig:springs_behavior}
  \end{subfigure}
  \caption{Periodic $1$D spring example schematic and constitutive
    behavior of the damage spring $0^{th}$ spring, a. Schematic of 1-D damage
    spring example, b. Force-displacement response of the damage spring}
  \label{fig:example_geom}
\end{figure}

As a simple mechanical system with non-convex energy functional, a periodic
$1$D spring system (schematic shown in~\autoref{fig:spring_schematic})
consisting of three nodes connected with springs ($k_0$, $k_1$, and $k_2$) is
taken as the first example.
The springs $k_1$ and $k_2$ are elastic springs with $k_1=k_2=k$,
while $k_0$ is a bi-linear damage spring, i.e. after a certain deformation
threshold $\gamma_0$ its mechanical behavior switches from elastic to
strain-softening.
The force-displacement response of the $k_0$ spring is depicted
in~\autoref{fig:springs_behavior}, The tangent of the strain-softening
phase of the constitutive behavior is $\alpha k$. Therefore,
the post-peak stiffness matrix of the system becomes:

\begin{equation}
  \label{eq:K_mat}
  \boldsymbol{K} =
  \begin{bmatrix}
    (1+\alpha)k & -k & -\alpha k \\
    -k     & 2k & -k \\
    -\alpha k   & -k & (1+\alpha)k
  \end{bmatrix}
\end{equation}
whose eigenvalues are $\left[ \lambda_1=0,\ \lambda_2=3k,\ \lambda_3= (2\alpha +1)k \right]$.
The third eigenvalue can be either positive or negative and
for values $\alpha < -\left(\rfrac{1}{2}\right)$, the system is not
\gls{psd} anymore, since it has one negative eigenvalue.

\begin{table}[]
  \centering
  \caption{Solvers used to solve 1D spring example}
  \label{tab:solver_springs_example}
  \begin{tabular}{|cc|ccc|}
    \hline
    \multicolumn{2}{|c|}{\multirow{2}{*}{Solver}}                                                                  & \multicolumn{3}{c|}{Functions needed as input}                           \\ \cline{3-5} 
    \multicolumn{2}{|c|}{}                                                                                         & \multicolumn{1}{c|}{\begin{tabular}[c]{@{}c@{}}Objective\\ (Energy)\end{tabular}} & \multicolumn{1}{c|}{\begin{tabular}[c]{@{}c@{}}Gradient\\ (Force)\end{tabular}} & \begin{tabular}[c]{@{}c@{}}Hessian\\ (Stiffness)\end{tabular} \\ \hline
    \multicolumn{1}{|c|}{\rom{1}} & \begin{tabular}[c]{@{}c@{}}Newton-CG\\ SciPy\end{tabular}                     & \multicolumn{1}{c|}{\xmark}     & \multicolumn{1}{c|}{\cmark}    & \cmark   \\ \hline
    \multicolumn{1}{|c|}{\rom{2}} & \begin{tabular}[c]{@{}c@{}}Trust Region \\ Newton-CG\\ Scipy\end{tabular}     & \multicolumn{1}{c|}{\cmark}     & \multicolumn{1}{c|}{\cmark}    & \cmark   \\ \hline
    \multicolumn{1}{|c|}{\rom{3}} & \begin{tabular}[c]{@{}c@{}}Modified\\  Trust Region\\  Newton-CG\end{tabular} & \multicolumn{1}{c|}{\xmark}     & \multicolumn{1}{c|}{\cmark}    & \cmark   \\ \hline
  \end{tabular}
\end{table}

The problem has been solved with
$k = 1.0,\ \gamma_0=0.1$ and for different values of $\alpha$. The boundary condition of
mean stretch equal to  $\overline{x}=0.11$, large enough to invoke post peak
behavior of the $k_0$ spring, is imposed. Three different solvers
listed in~\autoref{tab:solver_springs_example} are employed to solve the
equilibrium of the $1$D spring system. The functions needed to be
explicitly evaluated in these solvers' algorithms are noted
in~\autoref{tab:solver_springs_example}.

The main difference between the modified and the standard trust region solver,
as also noted in~\autoref{tab:solver_springs_example}, is that explicit
evaluation of the objective function is not needed in the modified solver.

The strain energy functional of the system as a function of the independent variable
$x_0$ for the imposed boundary condition of $\overline{x}=0.11$
is represented in~\autoref{fig:spring-response} for three different
values of $\alpha$, respectively from left to right, corresponding to convex,
meta-stable, and concave energy surfaces.
The variation of energy, and the final solution of the solvers
listed in~\autoref{tab:solver_springs_example} are depicted in this
figure. It is clear that in non-convex cases, Newton-CG solver is not capable
of finding the energy functional minimum, while both standard
and modified trust-region Newton-CG solvers converged to
minimizer of the energy (equilibrium points).
The hyper-parameters of the standard and the modified trust region solvers
(such as initial and maximum trust-region radius) are chosen to be identical,
this causes the solution steps of the solvers to coincide solving the $1$D
spring system.

The transition of the energy functional from convexity to non-convexity
is depicted in~\autoref{fig:spring-response} for $\alpha < -0.5$.
In~\autoref{fig:spring-response}a, the energy functional is
convex over all values of $x_0$ while in~\autoref{fig:spring-response}b
and ~\autoref{fig:spring-response}c,
the energy functional around the transition point $(x_0=0.1)$ of the spring
$k_0$ is non-convex.

This simple example can clearly show the equivalence of the obtained results with
that of the standard trust region algorithm. The availability of the energy
functional of this example makes the standard trust-region solver
applicable.
However, in general non-convex homogenization problems, the
energy functional is not always available and therefore the standard
trust-region solver is not an option and one can use the modified version
with the approximated energy functional.

{\centering
  \begin{figure}%[htbp]
    % \captionsetup{justification=centering}
    % \captionsetup{justification=centering,margin=0.1cm}
    \centering
    \includegraphics[width=12.5cm]{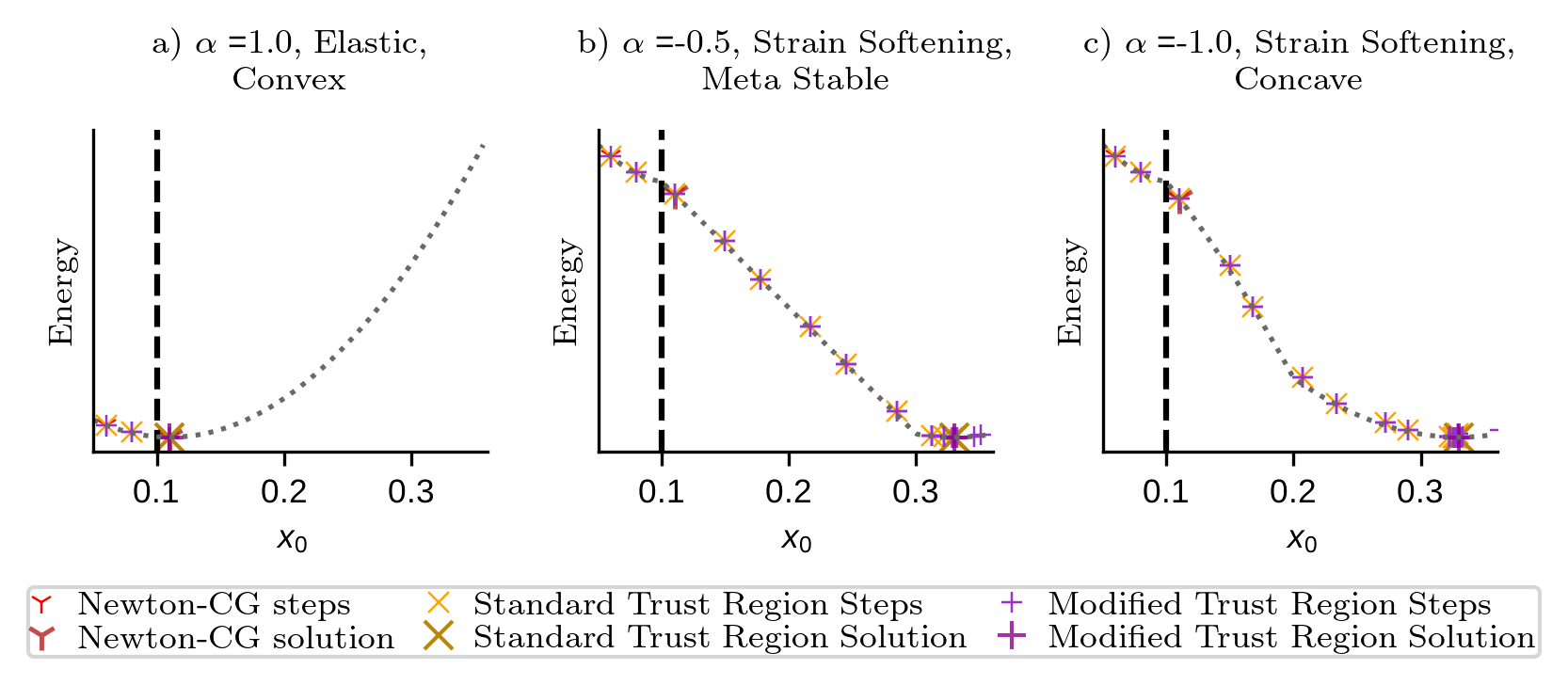}
    \label{fig:spring-response-energy}
    \caption{The energy of 1D spring system,
      schematically depicted in~\autoref{fig:spring_schematic}, over the solution
      trajectory of different solvers as a function of  $x_0$, for (a) Convex
      $\alpha=1$ (b) Meta-stable $\alpha=-\rfrac{1}{2}$, (c) non-convex $\alpha=-1$
      problems. In non-convex cases ((b), (c)), there is a concave point
      located at $x_0=0.1$ as the strain-softening behavior of the $k_0$
      springs is activated from that point on. Dashed line is the energy
      functional of the system as a function of $x_0$}
    \label{fig:spring-response}
  \end{figure}
}

\subsection{Convex example}
\label{sec:fig:convex_example_eshelby}
In order to examine the introduced modified trust-region solver solving convex
problems, a small-strain Eshelby inhomogeneity elasticity problem is chosen here
as the second numerical example. The Eshelby inhomogeneity is an
ellipsoidal body embedded in an infinite elastic medium, where the material
properties of the inhomogeneity differ from those of the matrix. The analytical
solution of the problem is known~\cite{Eshelby1957, Eshelby1959, Mura1982, Meng2012}.
A 2D example identical to the Eshelby inhomogeneity example presented
by~\citet{Leute2022} is considered here as our second example.
The linear \gls{fe} discretization of \citet{Leute2022} is adopted and
the problem is solved by two solvers, namely the plain Newton-\gls{cg} solver and
our modified trust region Newton-\gls{cg} solver.

\begin{figure}
  \centering
  \includegraphics[]{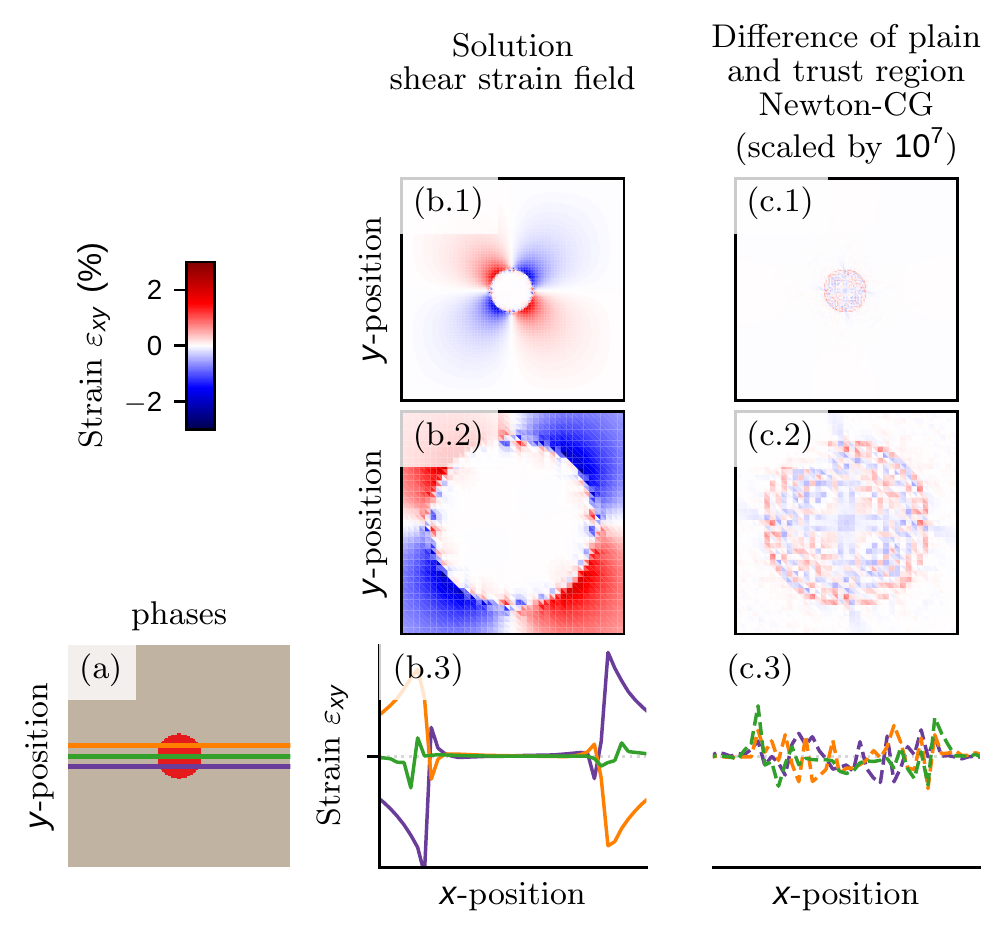}
  \label{example_geom:2D}
  \caption{Solution of the cylindrical Eshelby inhomogeneity
    problem under mean spherical strain of $0.01$ ($\varepsilon_{xx}=\varepsilon_{yy} =
    0.01, \varepsilon_{xy}=0$). Column a. shows the phase distribution of a soft
    inhomogeneity cylindrical Eshelby inhomogeneity problem (showing the inhomogeneity
    in red). Column b. shows solution of cylindrical Eshelby
    inhomogeneity with Newton-CG and Trust region Newton-CG as they look the same.
    The column (c) consists the difference of Newton-CG and Trust region
    Newton-CG solutions scaled by a factor of $10^7$ to make the difference
    visible. The first row shows the variation of the
    shear strain all over the solution domain. Second row shows the same
    variable zoomed around the inhomogeneity.  The third row shows
    the variation of shear strain over the green, purple and orange cuts
    (located at the center-line, $\rfrac{r}{2}$ below and above of the
    center-line of the inhomogeneity) in subfigure a. The third row also
    corresponds to the zoomed area around the inhomogeneity.}
  \label{fig:convex_example_eshelby}
\end{figure}

The schematic of the \gls{rve} geometry is presented in column (a) of
~\autoref{fig:convex_example_eshelby}.~\autoref{fig:convex_example_eshelby}
illustrates  the solution of both Newton-\gls{cg} and modified trust-region
Newton-\gls{cg} (column (b)). Column (c) consists the difference of the
solution of these two solvers.~\autoref{fig:convex_example_eshelby} depicts
that the solution of the two solvers are identical with respect to the solution
tolerance ($\eta_{NR}$ in~\autoref{algo:Trust region Newton-CG}) taken for the
iterative solution termination and the slight difference is in order of
magnitude of the tolerance.

\begin{figure}
  \centering
  \includegraphics[height=12cm]{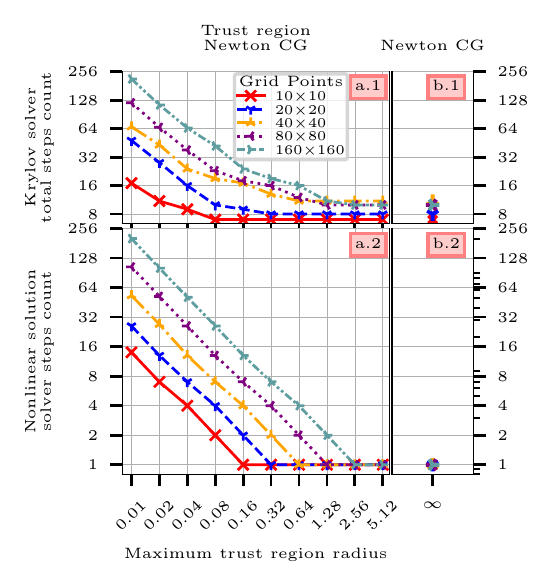}
  \label{example_geom:2D_steps}
  \caption{Number of Krylov steps (the first row) and number of nonlinear
    solution steps (the second row) needed for solving the Eshelby inhomogeneity
    problems for different number of grid points for the Newton-CG
    (column (a)) and as a function of initial trust region radius of trust region
    Newton-CG (column (b)). Number of nonlinear steps includes Newton steps and
    failed trial or trust region steps during the equilibrium solution}
  \label{fig:step_count_different_trust_radia}
\end{figure}

The number of Krylov solver and nonlinear solution steps needed to solve for
equilibrium versus trust radius variation (maximum trust region radius) is
plotted in~\autoref{fig:step_count_different_trust_radia} which shows that in a
convex problem the number of nonlinear solution steps as well as the
accumulative number of the Krylov solver steps needed to reach the solution
decays to that of Newton-CG solver as the size of the trust region
increases for all grid point counts.

\autoref{fig:step_count_different_trust_radia} also shows that, in
order to maintain the same number of nonlinear solution steps for solving a
problem with twice as many grid points in each direction (4 times
discretization points), the trust region should be roughly doubled. This
correlation is rooted in the fact that the trust region radius is actually the
radius of the hyper-sphere in the space of problem unknowns (strain in case of
projection-based solver). Imagine the discretization of a problem is refined by
a factor of $N$ in each spatial direction in a $2$D problem. This results in
$N^2$ scaling of the number of the discretization points. Accordingly, the size
of an equivalent solution step scales by $N$, in other words, an equivalent
step in the problem with refined discretization is $N$ times larger. As a
result, to maintain the ratio of the solution step length constant with respect
to trust region radius, The trust region radius should be scaled by a factor of
$N$. This finding suggests that in order to use the trust region solver
effectively one should loosen the trust region for larger problem
sizes. However, it should be noted that the trust region radius is actively
corrected during the solution of non-convex problems according to the accuracy
of the sub-problem model functional.

\subsection{Non-convex damage example}
\label{sec:non_convex_example}
As the third example, we will discuss an interesting real-world
example, namely~\gls{asr} damage homogenization.
\gls{asr} is one of the most widespread causes of internal concrete
deterioration~\cite{Hobbs1988, Swamy1991}. It is of great importance due to its
economical significance since it concerns critical structures
such as dams~\cite{Sellier2017}. \gls{asr} is initiated by a chemical reaction
between the alkali content of the cement paste and the silica in the
aggregates. Humidity is critical for initiation of this chemical reaction as it
washes the alkali content of the cement paste into the aggregate phase. The
resultant of this process (known as \gls{asr} gel) as a result of its
hydrophilic nature absorbs a considerable amount of water and therefore expands.
Due to their expansion, the confined \gls{asr} gel pockets subject their
surrounding to highly localized stresses.
The induced stress can damage the concrete micro-structure by creating microcracks
which will grow as the \gls{asr} advances. Following coalescing of
the cracks, advancement of \gls{asr} damage can result in macroscopic cracks
and stiffness and strength loss.

This example deals with \gls{rve} scale \gls{asr} damage modeling, in which we solve
equilibrium equation of a problem with a non-convex energy functional, where the
explicit expression of the functional is missing.
The extremely high computational costs of this problem, namely
\gls{asr} damage simulation has made it really challenging to conduct a
comprehensive multi-scale modeling of \gls{asr}~\cite{CubaRamos2018,
  Ramos2017}. In this example, we will present how the
modified trust region solver developed in this paper applied in the
projection-based scheme with~\gls{fe} discretization can offer a robust non-convex
\gls{rve} solver which suites \gls{asr} damage simulation and with its
quasi-linear scaling can make efficient multi-scale \gls{asr} simulation
possible.

\begin{figure}
  \captionsetup{}
  \centering
  \begin{subfigure}{.5\textwidth}
    \caption{}
    \captionsetup{justification=centering}
    \centering
    % \begin{overpic}
    %   \put(20, 100){A}
    % \end{overpic}

    \includegraphics[height=6.0cm, width =
    6.0cm]{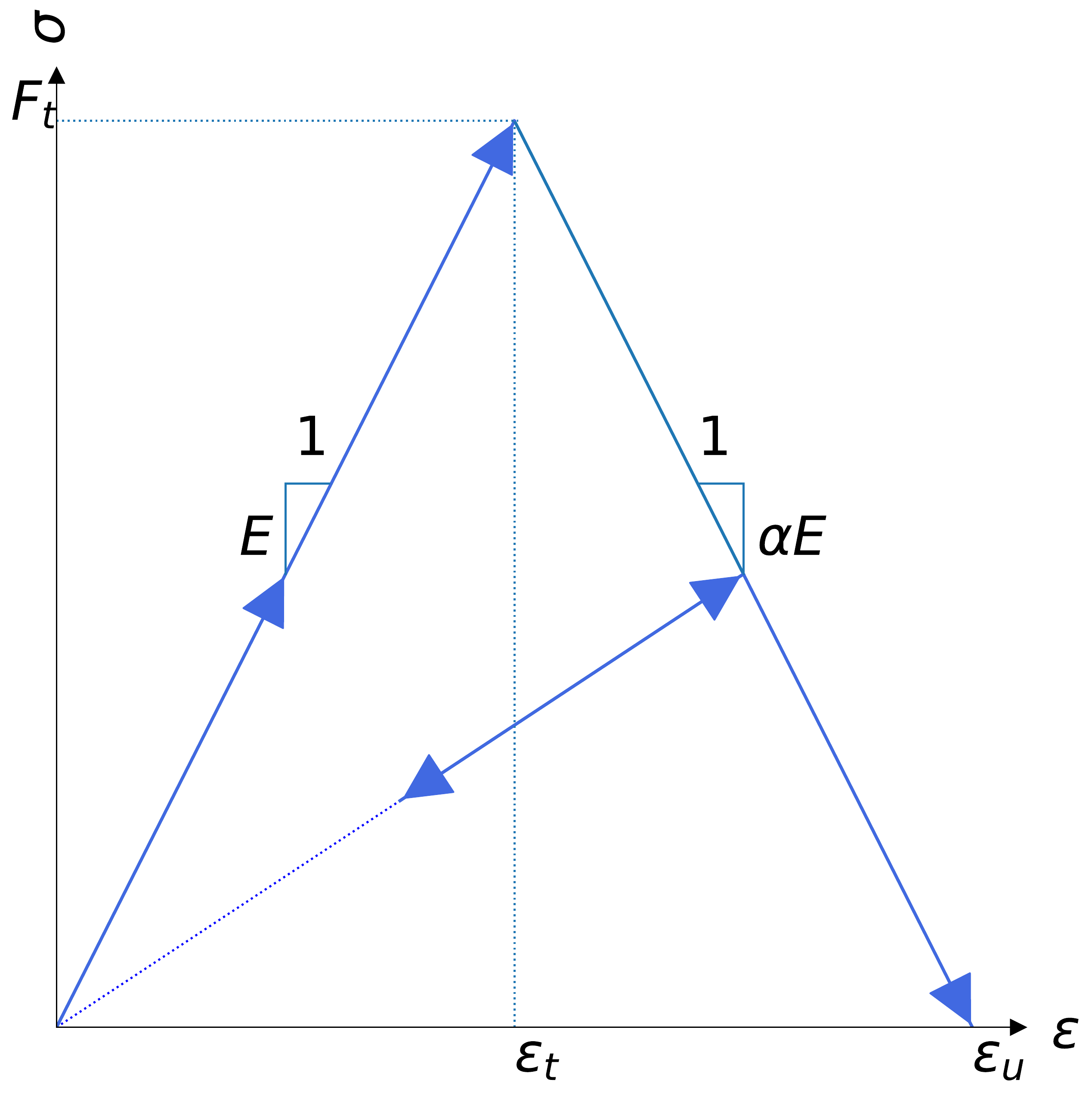}
    \label{fig:strain_stress}
  \end{subfigure}%
  \begin{subfigure}{.5\textwidth}
    % \captionsetup{justification=centering}
    \caption{}
    \centering
    % \begin{overpic}
    %   \put(20, 100){B}
    % \end{overpic}
    \includegraphics[height=6.0cm, width = 6.0cm]{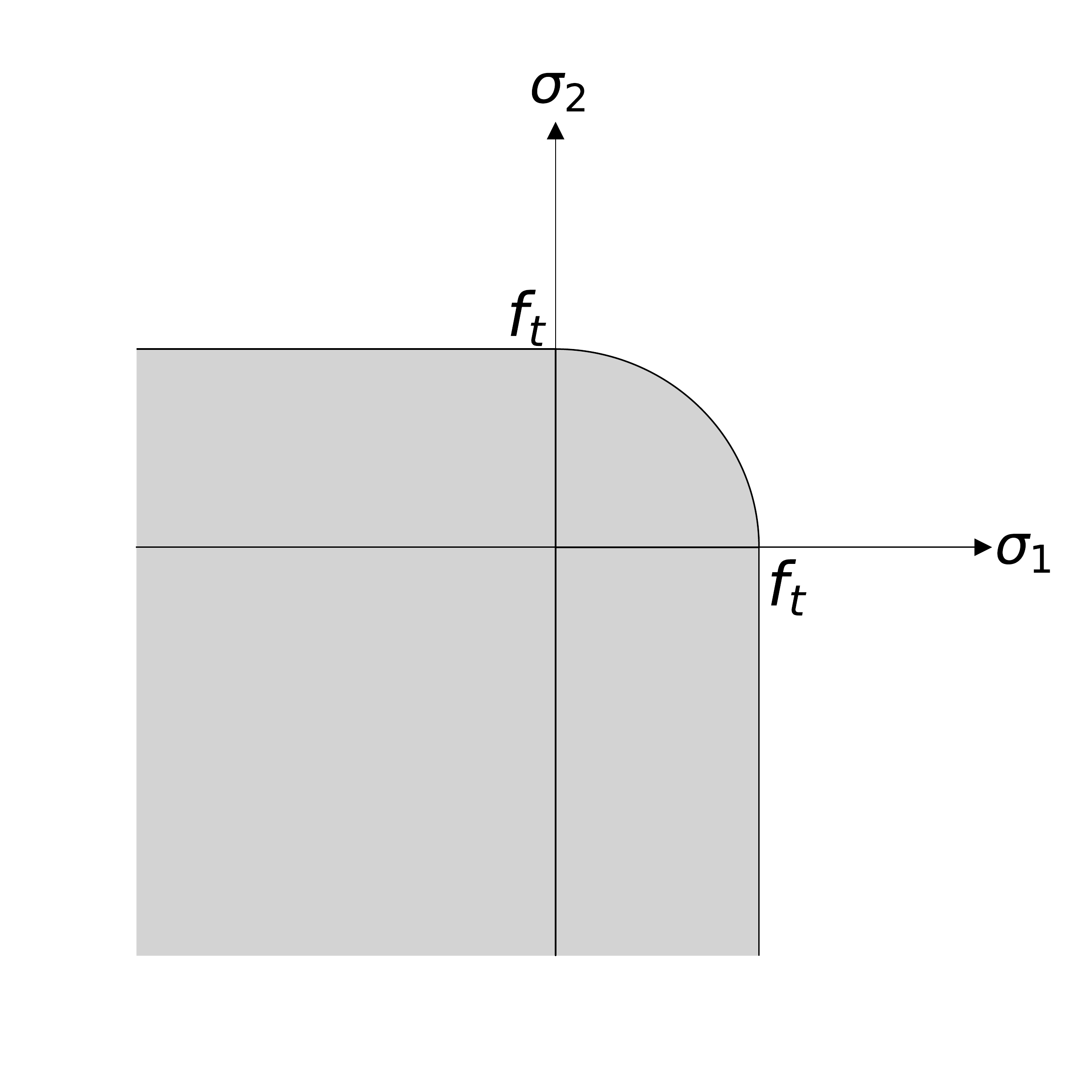}
    \label{fig:failure_criterion}
  \end{subfigure}
  \caption{Constitutive law of the damage material in ASR damage example,
    a. stress-strain response of the damage material and b.Failure criterion of the damage material
  }
  \label{fig:constitutive_law}
\end{figure}

In order to be able to compare our results with available \gls{asr} damage
meso-scale simulation results~\cite{Ramos2017}, a $2$D meso-scale \gls{asr}
damage model \gls{rve}s is considered here. We also report
the \gls{rve} stiffness reduction, as an important measure for expressing the
extent of \gls{asr} damage development~\cite{Ramos2017}.
In order to have same boundary conditions to~\cite{Ramos2017}
we imposed mean stress value as the boundary condition of
the \gls{rve}. In order to impose mean stress value, we applied the
necessary changes on the projection operator $\mathbb{G}$ according
to~\cite{Lucarini2019}.

\begin{figure}
  \includegraphics[height=16.0cm, width =
  12.0cm]{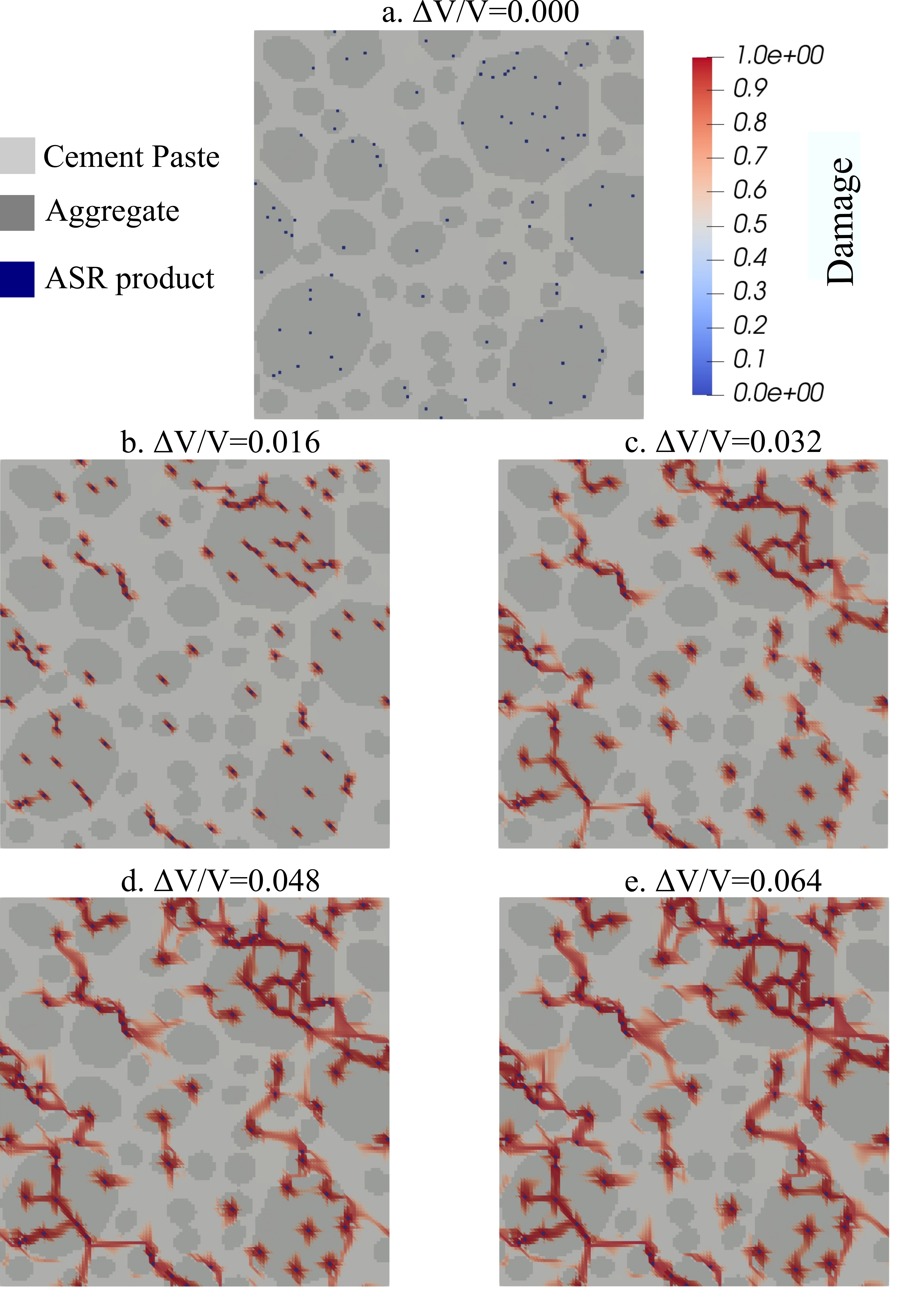}
  \caption{Evolution of crack pattern in the concrete micro-structure as a
    result of application of eigen-strain at ASR product sites shown as dark
    blue dots in the subfigure a. subfigures from a to e are arranged in a
    chronological sense and show snapshots of crack evolution inside the
    micro-structure as a function of increasing
    eigen-strain.}
  \label{fig:crack_pattern}
\end{figure}
In $2$D \gls{rve}s of this example, cement paste and aggregates are
explicitly resolved as two different phases.
Aggregates are placed inside a cement paste matrix in \gls{rve}s according to
Fuller size distribution~\cite{MotahariTabari2018}.
Pixels considered to be  containing growing \gls{asr} gel pockets has been
randomly inserted inside aggregates. The structure of the $2$D micro-structure
is depicted in~\autoref{fig:crack_pattern}a.

The constitutive laws of both aggregate and
cement paste phases are a bilinear crack band damage with an isotropic damage
measure~\cite{Mazars1989}.
The initial Young modulus of $E^0$ and its damaged counterpart ($E$) can be
related through:
\begin{equation}
  \label{eq:energy_damage_dunant}
  E = (1-D)E^0,
\end{equation}
where the damage variable
$D$ can vary between $0$ for intact material to $1$ for completely damaged
(fractured) material. Further details of the damage constitutive law
is presented in~\ref{app:stiffness}.
This constitutive behavior is also depicted
schematically in~\autoref{fig:constitutive_law} where the subfigures (a) and
(b) respectively illustrate the stress-strain response and the failure criteria
of the constitutive law. As illustrated in~\autoref{fig:failure_criterion}, the
damage material only fails under tensile loads.

\begin{table}[]
  \centering
  \caption{Material properties of the non-convex damage example obtained from
    Ref.~\cite{Gallyamov2020} for characteristic size $l_c = 5 \times 10^{-4}m$}
  \label{tab:mat_props}
  \begin{tabular}{lccccc}
    \hline \\ [-0.7em]
    & E {[}GPa{]} & $\mu${[}GPa{]} & $\nu$ & $G_{\mathrm{c}}${[}J/$m^2${]} & $f_t^0$
                                                              {[}MPa{]}
    \\ \\[-0.7em] \hline\hline
    Aggregates   & 59          & 22.6           & 0.3   & 160                & 10                \\
    Cement paste & 12          & 4.6            & 0.3   & 60                 & 3                 \\
    ASR product  & 11          & 4.7            & 0.18  & -                  & -                 \\ \hline
  \end{tabular}
\end{table}

The parameters of the constitutive law, listed in~\autoref{tab:mat_props}
are obtained from Ref.~\cite{Gallyamov2020}.
The damage part of this constitutive law leads to a \gls{snpsd} system
matrix. For different discretization sizes, we conduct energy-based
constitutive law regularization to preserve the Mode-I fracture energy
($G_c$) of the damage constitutive law and make the solution of the equilibrium
independent of the mesh size.

We model the expansion of the \gls{asr} product sites by applying
eigenstrain on the pixels containing them. Growing eigenstrain is added to the
strain associated to these quadrature points before their constitutive law
evaluation. These specific pixels are modeled as a linear elastic phase
and assumed to contain the growing \gls{asr} gel pockets inside them.
Line $28$ of the~\autoref{algo:Trust region Newton-CG} is where eigenstrain is
actually applied on the system. The mean value of stress is imposed as a
boundary condition to be zero (free expansion) for the \gls{rve}. Imposing mean
stress value as the boundary condition is implemented by means of modifying
the projection operator following~\cite{Lucarini2019}.
This choice of boundary condition allows us to compare the results to similar
results in the literature.

\begin{figure}
  \includegraphics[]{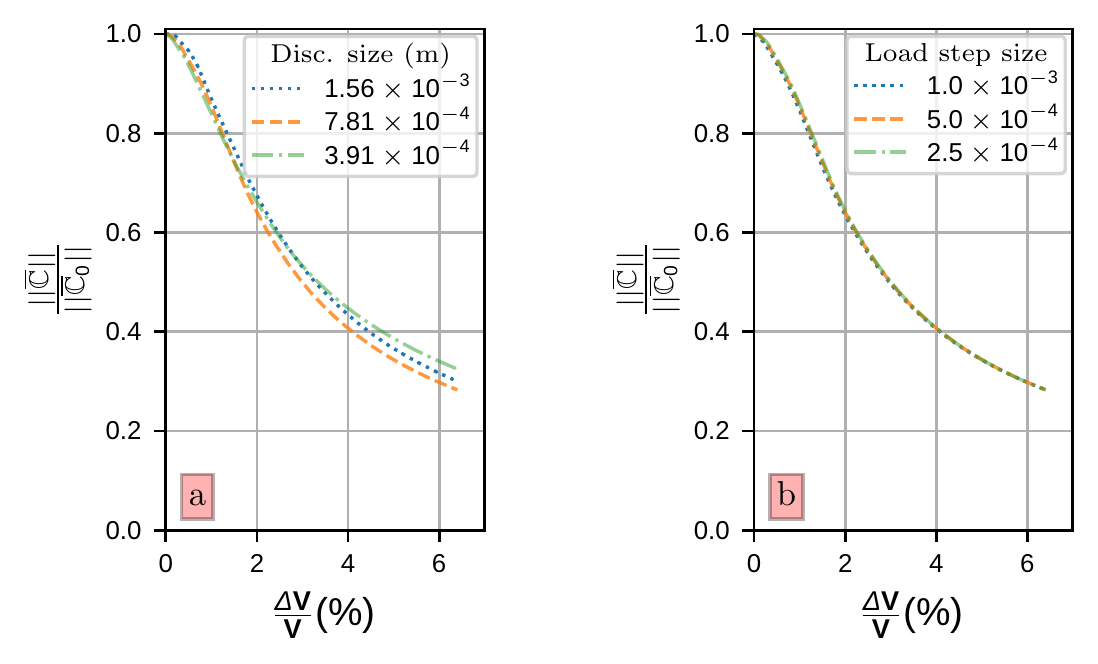}
  \caption{Discretization (mesh size) and load step study of the stiffness
    deterioration of $2$D ASR damage simulation expressed as the ratio of the
    norm of the effective stiffness $||\overline{C}||$ tensor during~\gls{asr}
    advancement divided by the norm of the effective stiffness tensor of the
    intact~\gls{rve} $||\overline{C}_0||$. in subfigure a) the stiffness
    reduction of the same problem with fine ($h_{\mathrm{f}}=3.91\times10^{-4}m$),
    medium ($h_{\mathrm{m}}=\ 2h_{\mathrm{f}}$), and coarse ($h_{\mathrm{c}}=\
    4h_{\mathrm{f}}$) grid carried out with the medium load step size ($\Delta
    \varepsilon_m^{\mathrm{eig}}= 5.0\times10^{-4}$) is plotted. In subfigure b) the stiffness
    reduction of the same problem with Load step study with small ($\Delta
    \varepsilon_{\mathrm{s}}^{\mathrm{eig}} = 2.50\times10^{-4}$), medium ($\Delta
    \varepsilon_{\mathrm{m}}^{\mathrm{eig}} =2 \Delta \varepsilon_{\mathrm{s}}^{\mathrm{eig}}$) and
    large load step of ($\Delta \varepsilon_{\mathrm{l}}^{\mathrm{eig}} =4 \Delta
    \varepsilon_{\mathrm{s}}^{\mathrm{eig}}$) carried out on $h_{\mathrm{m}}$. is
    plotted.
  }
  \label{fig:size_study}
\end{figure}

Crack pattern advancement inside a sample \gls{rve} of physical size of
$0.1 \mathrm{m}\times 0.1 \mathrm{m}$ is depicted in~\autoref{fig:crack_pattern}.
As shown in~\autoref{fig:crack_pattern}b, the damage initiation sites are adjacent
to the growing gel pixel sites in the \gls{rve}. The advancement of the cracks
caused by~\gls{asr} damage is depicted in~\autoref{fig:crack_pattern}c-e.
These subfigures illustrate the crack coalescence process as the~\gls{asr}
expansion proceeds. Crack coalescence depends on the
distance between the gel pocket and the aggregate boundary as well as the
distance to other gel pockets. Crack coalescence occurs earlier at crack sites
near the boundary of the aggregates or other gel pockets.

% An analytical approach
% derived from Hill-Mandel lemma~\cite{Nicot2017} and inspired by Chapter 4 of
% ~\cite{Yvonnet2019} is developed and utilized here to obtain the
% effective stiffness of the \gls{rve}. In this approach,
% the local stiffness of the discretization points at solved equilibrium is
% utilized to calculate the effective stiffness without need of applying and
% solving for compression or tensile test load steps.
Next, we conducted the discretization study and the load step study on a
\gls{rve} with free expansion and growing gel pixels. In order to be able to
realize nominally equivalent load scenarios by variation of the discretization
size we have kept the area expansion of the pixels containing \gls{asr} gel
pockets constant.

The effective stiffness reduction for a micro-structure subjected
to similar area expansion of gel pocket pixels is presented with different
discretization in~\autoref{fig:size_study}a and load step size
in~\autoref{fig:size_study}b.
~\autoref{fig:size_study}a shows that the results
obtained are independent of discretization size as the result of
fine, medium and coarse grain are matching. The stiffness reduction of the
medium grid size with different load step size application illustrates that
the results are also independent of load step size.

\begin{figure}
  \centering
  \centering \includegraphics[width=12cm]{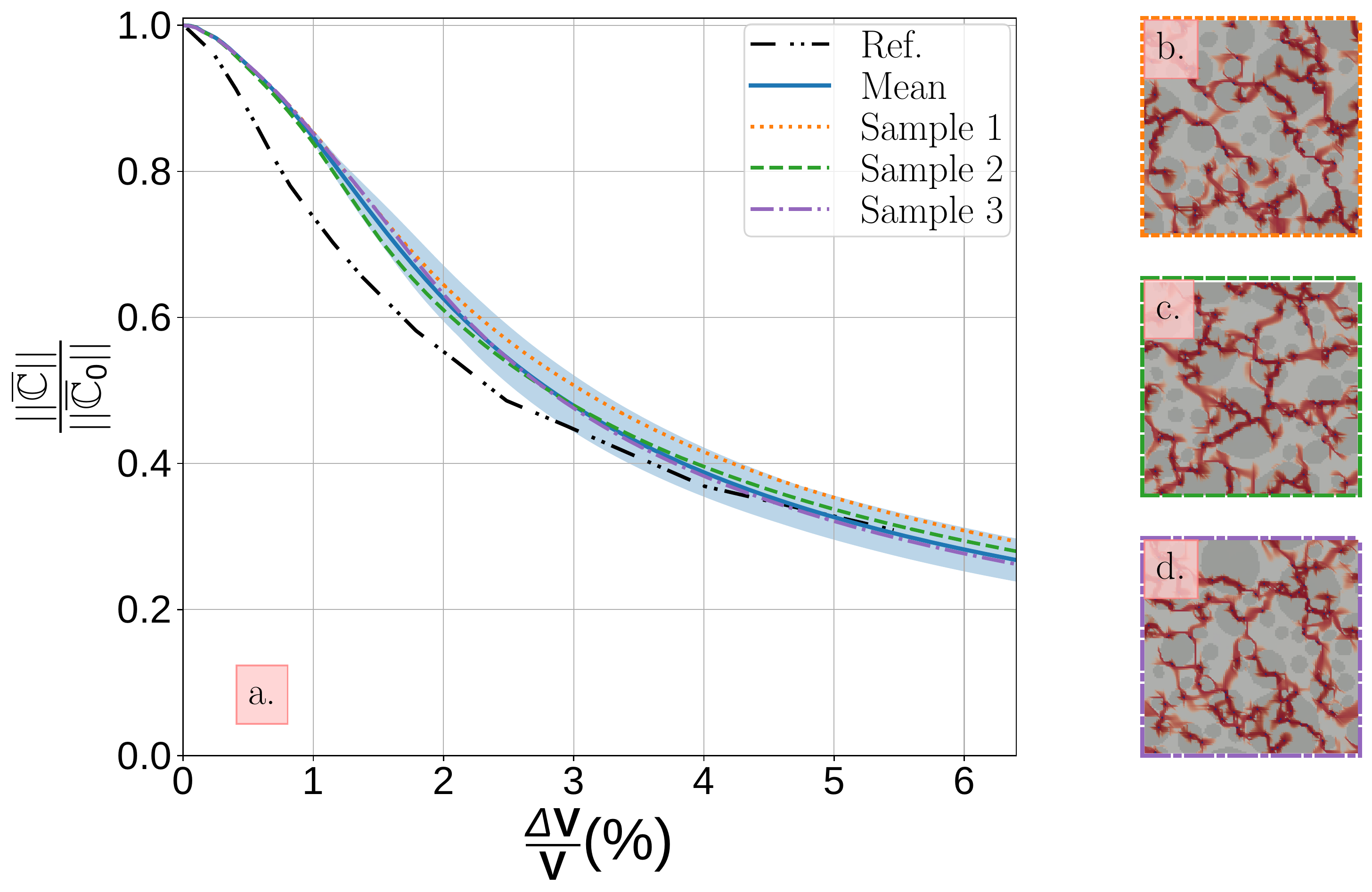}
  \caption{Stiffness reduction of sample micro-structrues are shown in subfigure
    a). The light blue area
    is the area shaped between the maximum and minimum stiffness reduction of
    100 samples. The solid blue line is the mean value of the stiffness
    reduction and the other 3 lines are the results corresponding to three
    represetative micro-structures. Subfigures b. to d. show the final crack
    pattern developed in the micro-structure corresponding to the
    three samples (sample 1, 2 , and 3) noted in the subfigure a.}
  \label{fig:differ_random_phases}
\end{figure}

According to the results of the discretization study we have chosen
the medium grid size and a load step size of $5.00 \times 10^{-4}$ (the medium value
of the considered load steps) to conduct the simulations with different
randomly generated micro-structures in the following of the paper whose results
are depicted in~\autoref{fig:differ_random_phases}. We then subjected
$100$  randomly generated concrete $2$D micro-structures to
\gls{asr} expansion (imposed as eigenstrain in pixels containing gel pockets)
under free expansion boundary condition.
The shaded blue area in~\autoref{fig:differ_random_phases}a shows the
distribution of the stiffness loss of these micro-structure vs imposed
volumetric expansion of the gel pixels. The mean of the distribution is also
plotted as solid blue line. The stiffness loss of $3$ representative
micro-structures are plotted in~\autoref{fig:differ_random_phases}a and the
\autoref{fig:differ_random_phases}b-d are their corresponding final crack
pattern.

In addition, the results obtained by~\citet{Ramos2017} modeling \gls{asr}
damage in a similar configuration though using \gls{fem} scheme and using a
different approach for addressing the numerical instabilities due to the
non-convexity of the problem (namely \gls{sla}) is also plotted
in~\autoref{fig:differ_random_phases}a labeled as Ref which depicts reasonable
agreement with our obtained results. By using our modified trust region solver
in the projection-based scheme with~\gls{fe} discretization
our solution is much faster than their approach.
According to correspondence with the authors of~\cite{Ramos2017,
  Gallyamov2020} their calculations on 28 cores take roughly 48 hours while our
simulations, on average, take half an hour on 16 cores on the same machine
(Fidis cluster at \gls{epfl}) which shows significant improvement on a $2$D
\gls{asr} damage simulation.
The small differences visible in~\autoref{fig:differ_random_phases} is probably due to
subtle distinctions between the models, namely we have used rectangular
elements containing \gls{asr} gel pockets while they were triangular in their
model. The aggregates in their model are assumed to have circular geometry
while they were ellipsoidal in our model. The differences in the application of
the boundary conditions can also be source of difference between the
models as they have traction free boundary conditions while we have applied
zero mean stress on our~\gls{rve}. However, despite of  all this subtle
differences our obtained stiffness loss is in a good agreement with their results.

\section{Summary and conclusion}

In this paper, we present a first order incremental approximation of the strain
energy functional, which makes it possible to use fast ringing-free
spectral solvers for non-convex problems, such as damage mechanics. We have
derived the approximated energy functional based on the Taylor expansion of the energy
functional of the system. Using the introduced incremental energy functional
enables employing modern optimization techniques such as quasi-Newton and
trust-region solvers in problems without easy access to the explicit objective function, for
instance non-linear mechanical homogenization problems. As an example,
we show, in this paper, how the introduced incremental energy functional makes the
use of trust-region Newton-\gls{cg} for computational homogenization possible.

We used the modified trust region sovler to solve a minimal $1$D non-convex
problem where the obtained solution was shown to be in agreement with that of
standard trust region Newton-\gls{cg} solvers. Next, we used it to solve a
convex homogenization problem; finally, we simulated and studied a real-world
homogenization problem with non-convex energy functional, namely meso-scale \gls{asr} damage,
was simulated and studied by means of the presented modified trust region
solver. The obtained results were compared to those reported in the literature~\cite{Ramos2017}
with much lower computational cost which made our solution roughly $200$ times faster.

\newpage
\newpage
%%%%%%%%%%%%%%%%%%%%%%%%%%%%%%%%%%%%%%%%%%%%%%%%%%%%%%%%%%%%%%%%%%%%%%%%%%%%%% 
\begin{algorithm}[H]
  \caption{Pseudo-algorithm of CG-Steihaug solver with reset}
  \label{algo:CG-Steihaug solver}
  \begin{algorithmic}[1]
    \State{\textbf{Solve for $\vek{r}$ with system Matrix $\tangent$ and initial RHS $\vek{b}_0$}}
    \State{$\eta_\text{CG}$ }\Comment{\gls{cg} tol.}
    \State{$j_\text{CG,max}$} \Comment{max iterations of \gls{cg}}
    \State{Set $\vek{r}_0\leftarrow \vek{0}, \vek{b}_0\leftarrow\vek{b}, \vek{d}_0\leftarrow-\vek{b}_0$}
    \Comment{initialization}
  \item[]
    \If{$\norm{\vek{b}} \leq \eta_\text{cg}$}
    \State{Return $\vek{p}_i=\vek{r}_0=0$} \Comment{already at solution}
    \EndIf{}
  \item[]
    \For{$j=0,1,2, ...,j_\text{CG,max}$}

    \If{$\vek{d}_j^T \tangent_i \vek{d}_j \leq 0$} \Comment{non-convex}
    \State{find $\tau$ such that $\vek{p}_i = \vek{r}_j + \tau \vek{d}_j$} minimizes
    $m_i(p_i)$
    \Statex{\qquad \qquad and satisfies $\norm{\vek{p}_i} = R_i$}
    \State{Return $\vek{p}_i$}
    \EndIf{}
  \item[]
    \State{$\alpha_j \leftarrow\vek{r}_j^T\vek{r}_j / \vek{d}_j^T \tangent_i \vek{d}_j$}
    \State{$\vek{r}_{j+1}\leftarrow\vek{r}_{j} + \alpha_j \vek{d_j}$}\Comment{update the iterate}

  \item[]
    \If{$\norm{\vek{r}_{j+1}} \geq R_i$} \Comment{hit the boundary of trust region}
    \State{find $\tau$ such that $\vek{p}_i = \vek{r}_j + \tau \vek{d}_j$} minimizes
    $m_i(p_i)$
    \Statex{\qquad \qquad and satisfies $\norm{\vek{p}_i} = R_i$}
    \State{Return $\vek{p}_i$}
    \EndIf{}
  \item[]
    \If{$\norm{\vek{r_{j+1}}} \leq \eta_\text{cg}$} \Comment{convergence satisfied}
    \State{Return $\vek{p}_i=\vek{r}_{j+1}$}
    \EndIf{}
  \item[]
    \If{$\vek{r}_{j+1} \cdot \vek{r}_{j} / \vek{r}_{j+1} \cdot \vek{r}_{j+1} >
      0.2$}\Comment{successive steps are not conjugate}
    \State{$\vek{r}_{j+1} = \tangent\vek{r_j} - \vek{b}_j$} \Comment{reset CG}
    \State{$\beta_{j+1} \leftarrow$ 0} \Comment{reset CG}

    \Else
    \State{$\beta_{j+1} \leftarrow \vek{r}_{j+1}^T\cdot \vek{r}_{j+1}/\vek{r}_{j}^T\cdot \vek{r}_{j}$}
    \EndIf{}
  \item[]
    \State{$\vek{d}_{j+1} \leftarrow -\vek{r}_{j+1} + \beta_{j+1} \vek{d}_{j}$}
    \Comment{compute new update direction}
    \EndFor{}
  \end{algorithmic}
\end{algorithm}

%%%%%%%%%%%%%%%%%%%%%%%%%%%%%%%%%%%%%%%%%%%%%%%%%%%%%%%%%%%%%%%%%%%%% 

\newpage
\newpage
%%%%%%%%%%%%%%%%%%%%%%%%%%%%%%%%%%%%%%%%%%%%%%%%%%%%%%%%%%%%%%%%%%%%%%%%%%%%%% 
\begin{algorithm}[H]
  \caption{Pseudo-algorithm of projection-based Newton-CG trust-region
    algorithm in small strain formulation}
  \label{algo:Trust region Newton-CG}
  \begin{algorithmic}[1]
    \State{\textbf{Initialize: }}
    \State{$\eta_{\text{eq.}},\ \eta_{\text{NR}}, \ \eta_{\text{CG}}, \ \eta_{\text{up.}}$ }\Comment{equilibrium-, Newton, \gls{cg} and update tol.}
    \State{$i_\text{NR,max},\ i_\text{CG,max}$} \Comment{max iterations Newton-Raphson and \gls{cg}}
    \State{$R, R_{max}$} \Comment{trust region radius, maximum radius}
    \State{$\strain = \t{0}$}
    \Comment{small-strain initial guess}

  \item[]
    \For{$\Delta \strain = \Delta \strain_1, \Delta \strain_2, \dots $}  \Comment{macroscopic strain increments}
    \State{$\strain=\strain+\Delta \strain$} \Comment{increment grad with load step}
    \State{$\strain_{eval} = \strain + \strain_{\mathrm{eig}}$} \Comment{adding eigenstrain if needed}
    \State{$\stress,\ \tangent = \stress(\strain_{eval}),\ \tangent(\strain_{eval})$}\Comment{evaluate stress and tangent}
    \State{$\vek{b}=-\mathbb{G}:\stress(\strain_{eval})$} \Comment{RHS calculation}

    \If{$\norm{b} \leq \eta_\text{eq.}$}
    \State{$\ \ $\text{Newton-Raphson converged}}
    % \item[]
    \Statex \hspace{\algorithmicindent} \hspace{\algorithmicindent} {Go to line 6} \Comment{linear problem, next load step}
    \EndIf{}

    % \item[]
    \For{$i = 0, 1, 2, \dots, i_\text{NR,max} $} \Comment{Newton-Raphson iteration}
    \State{\text{Prepare coefficient matrix of the linearized equation $ \mathbb{G}:\tangent:\delta \strain = b$}}
    \State{\textbf{Solve $ \mathbb{G}:\tangent:\delta \strain = b$ for $\delta \strain$ with Steihaug \gls{cg}~\cite{Nocedal2006}:}}
    \Statex{$\qquad \qquad \qquad$\text{in} $i_\text{CG,max}$ \text{steps to accuracy}
      $\eta_\text{CG}$,~\autoref{algo:CG-Steihaug solver}}
    \State{$\Delta m_i=\stress:\delta \strain + 1/2\ \delta\strain:\tangent:\delta\strain$} \Comment{ energy model change}
    \State{$\stress^{\mathrm{trial}}=\stress(\strain+\delta \strain)$}\Comment{stress evaluation with trial strain}
    \State{$\overline{\Delta W} = 1/2\ (\stress:\delta \strain + \stress^{\mathrm{trial}}:\delta \strain)$} \Comment{$1^{st}$order energy approx. change}
    \State{$\overline{\rho}=\overline{\Delta W}/\Delta m_i$}
    \If{$\overline{\rho} < 1/4$}
    \State{$R \leftarrow 1/4\ R$}\Comment{shrink trust region}
    \Else
    \If{$\overline{\rho} > 3/4$ and $||\delta \strain|| = R$}
    \State{$R \leftarrow \min(2R,\ R_{max})$} \Comment{expand trust region if possible}
    \EndIf{}
    \EndIf{}
    \If{$\overline{\rho} > \eta_{up.}$}
    \State{$\strain \leftarrow \strain + \delta \strain$} \Comment{increment grad with accepted
      solution step}
    \State{$r_{\text{NR}} =  \norm{\delta \strain} / \norm{\strain}$} \Comment{calculating relative residual}
    \State{$\strain_{\mathrm{eval}} = \strain + \strain_{\mathrm{eig}}$} \Comment{adding eigen strain if needed}
    \State{$\stress,\ \tangent = \stress(\strain_{\mathrm{eval}}),\ \tangent(\strain_{\mathrm{eval}})$}\Comment{evaluate stress and tangent}
    \State{$b=-\mathbb{G}:\stress(\strain_{\mathrm{eval}})$} \Comment{RHS calculation with updated grad}
    \If{$||\delta \strain|| < R$}
    \If{$\norm{b} \leq \eta_\text{eq.} \textbf{ or }  r_{\text{NR}} \leq \eta_\text{NR}$}
    % \State{$\strain_0 = \strain_{i+1}$}
    \State{\text{Newton-Raphson is converged}}
    \Statex \hspace{\algorithmicindent} \hspace{\algorithmicindent}
    \hspace{\algorithmicindent} \hspace{\algorithmicindent} {Go to line 6} \Comment{next load step}
    \Else
    \Statex \hspace{\algorithmicindent} \hspace{\algorithmicindent}
    \hspace{\algorithmicindent}  \hspace{\algorithmicindent} {Go to line 14} \Comment{next Newton loop iteration}
    \EndIf{}
    \Else
    \Statex \hspace{\algorithmicindent} \hspace{\algorithmicindent}
    \hspace{\algorithmicindent} {Go to line 14} \Comment{next Newton loop iteration}
    \EndIf{}
    \Else
    \State{\text{Trial step rejected}}
    \Statex \hspace{\algorithmicindent}  \hspace{\algorithmicindent} {$\ \ $  Go to line 14} \Comment{next Newton loop iteration}
    \EndIf{}
    \EndFor{} %Newton
    \EndFor{} %Strain increment
  \end{algorithmic}
\end{algorithm}

%%%%%%%%%%%%%%%%%%%%%%%%%%%%%%%%%%%%%%%%%%%%%%%%%%%%%%%%%%%%%%%%%%%%% 

\section*{Acknowledgments}
We acknowledge funding by the the Swiss National Science
Foundation (Ambizione grant 174105~(TJ)), European Research Council
(StG-757343~(LP)), the Carl Zeiss Foundation (Research cluster "Interactive and
Programmable Materials - IPROM"~(LP)), the Deutsche Forschungsgemeinschaft (EXC
2193/1 - 390951807~(LP)), the Czech Science Foundation (projects
No.~20-14736S~(ML) and 19-26143X~(JZ)), and the European Regional Development
Fund (Centre of Advanced Applied Sciences – CAAS,
CZ.02.1.01/0.0/0.0/16\_019/0000778 (ML, IP)).
We also would like to acknowledge Prof. Jan Zeman, Prof. Lars Patweska and
Prof. Ivana Pultarová for their help and guidance in writing the paper.

\appendix
\section{Damage material}
\label{app:stiffness}
The scalar damage measure $D$ used in~\eqref{eq:energy_damage_dunant} is a function of
the highest control variable $\kappa$ that the material had experienced in the
solution history. Different
definitions of the control variable $\kappa$ result in different damage models.
Here, for instance, both phases are assumed to be only damaging in tension. To
this end the damage control measure is considered as:

\begin{equation}
  \label{eq:strain_measure}
  \kappa = \norm{\strain^{(t)}},\ \mathrm{with}\ \strain^{(t)}=
  \mathcal{H} \left( \epsilon_i \right) \vek{q}_i,
\end{equation}
where $\mathcal{H}$ is the Heaviside function, $\epsilon_i$ is the $i^{\mathrm{th}}$
eigenvalue of the strain tensor and $\vek{q_i}$ is its corresponding eigenvector.
This failure criterion suites \gls{asr} damage simulation when cracking under
tension (Mode-I fracture) is the effective damage mechanism. The damage
variable evolves according to the flow rule of
\begin{align}
  \label{eq:damage_evolution}
  \dot{D} = 
  \begin{cases}
     \  0 & \mathrm{if} \ f < 0\  \mathrm{or}\  f = 0\
    \ \mathrm{and}\  \dot{f}<0 ,\
    \mathrm{and} \ \dot{\phi} <0, \\
    \ \Phi \left(\t{\strain}, D \right) & \mathrm{if}\ f=0,
    \ \mathrm{and} \ \dot{f}=0 \\
  \end{cases}
\end{align}
where $\Phi  \left(\t{\strain}, D \right)$, called the damage evolution
function is a positive function of the strain and damage variable. The damage
loading surface $f$ corresponds to the damage threshold of the material.
In the bilinear scalar damage evolution considered here, the damage variable
can be calculated according to:
\begin{equation}
  \label{eq:our_model_damage_variable}
  D = \frac {\left(\kappa - \kappa_0\right)\left(1+\alpha\right)}{\kappa}
\end{equation}
where $\kappa_0$ is the initial strain measure damage threshold and $\alpha$ is
the relative slope of the softening part of the constitutive law, depicted
in~\autoref{fig:strain_stress}.

It is notable that computing the tangent stiffness of the
constitutive law, needed due to using Newton-\gls{cg} solvers, can introduce
singularity because it involves differentiation of $\kappa$ with respect to
$\strain$. According to~\eqref{eq:strain_measure} this differentiation
needs differentiation of strain eigenvectors ($\vek{q}_i$) with respect to the
strain tensor itself with the form of:
\begin{equation}
  \label{eq:eigenvector_diff}
  \frac{\partial \vek{q}_{i,\gamma}}{\partial \strain_{\alpha \beta}}= \sum_{i \neq j}{\frac{\vek{q}_{i, \beta}\vek{q}_{j, \alpha}}{\epsilon_i -
      \epsilon_j} \vek{q}_{j, \gamma}}
\end{equation}
which can cause singularity in case $\epsilon_i-\epsilon_j$ tends to $0$.
In order to resolve this problem, we have reformulated the $\strain^{(t)}$ by
defining a so called masking matrix $\t{M}^{(t)}$ based on the spectral
decomposition of $\strain$ (inspired by~\cite{Contrafatto2007}) according to:
\begin{align}
  \label{eq:stress_contributions_maksing_matrices}
  \strain^{(t)} =\ \t{M}^{(t)}\  \strain\  \ \t{M}^{(t)}
\end{align}
where  $\t{M}^{(t)}$ is defined as:
\begin{align}
  \label{eq:Mask_matrices_t}
  \t{M}^{(t)}  = \sum_{i=1}^{d}\ {\mathcal{H} \left( \epsilon_i \right) \vek{q}_i \otimes
  \vek{q}_i}.
\end{align}
Working out the differentiation of $\strain^{(t)}$ expression
with respect to $\strain$ according
to~\eqref{eq:stress_contributions_maksing_matrices} and
~\eqref{eq:Mask_matrices_t} singularity is
avoided in the explicit stiffness tangent formulation.

\bibliographystyle{abbrvnat}
\bibliography{main}

\begin{thebibliography}{81}
\providecommand{\natexlab}[1]{#1}
\providecommand{\url}[1]{\texttt{#1}}
\expandafter\ifx\csname urlstyle\endcsname\relax
  \providecommand{\doi}[1]{doi: #1}\else
  \providecommand{\doi}{doi: \begingroup \urlstyle{rm}\Url}\fi

\bibitem[Bažant(1976)]{Bazant1976}
Z.~P. Bažant.
\newblock Instability, ductility, and size effect in strain-softening concrete.
\newblock \emph{Journal of the Engineering Mechanics Division}, 102\penalty0
  (2):\penalty0 331--344, apr 1976.
\newblock \doi{10.1061/jmcea3.0002111}.

\bibitem[Budiansky(1965)]{Budiansky1965}
B.~Budiansky.
\newblock On the elastic moduli of some heterogeneous materials.
\newblock \emph{J. Mech. Phys. Solids}, 13\penalty0 (4):\penalty0 223--227,
  1965.

\bibitem[Byrd et~al.(2000)Byrd, Gilbert, and Nocedal]{Byrd2000}
R.~H. Byrd, J.~C. Gilbert, and J.~Nocedal.
\newblock A trust region method based on interior point techniques for
  nonlinear programming.
\newblock \emph{Math. Program.}, 89\penalty0 (1):\penalty0 149--185, 2000.

\bibitem[Conn et~al.(2000)Conn, Gould, and Toint]{Conn2000}
A.~R. Conn, N.~I. Gould, and P.~L. Toint.
\newblock \emph{Trust region methods}.
\newblock SIAM, 2000.

\bibitem[Contrafatto and Cuomo(2007)]{Contrafatto2007}
L.~Contrafatto and M.~Cuomo.
\newblock Comparison of two forms of strain decomposition in an elastic-plastic
  damaging model for concrete.
\newblock \emph{Modelling and Simulation in Materials Science and Engineering},
  15\penalty0 (4):\penalty0 S405--S423, May 2007.
\newblock ISSN 1361-651X.
\newblock URL \url{http://dx.doi.org/10.1088/0965-0393/15/4/S07}.

\bibitem[Cuba~Ramos(2017)]{Ramos2017}
A.~I. Cuba~Ramos.
\newblock \emph{Multi-Scale Modeling of the Alkali-Silica Reaction in
  Concrete}.
\newblock PhD thesis, EPFL, Lausanne, 2017.
\newblock URL \url{http://infoscience.epfl.ch/record/227479}.

\bibitem[Cuba~Ramos et~al.(2018)Cuba~Ramos, Roux-Langlois, Dunant, Corrado, and
  Molinari]{CubaRamos2018}
A.~I. Cuba~Ramos, C.~Roux-Langlois, C.~F. Dunant, M.~Corrado, and J.-F.
  Molinari.
\newblock {HPC} simulations of alkali-silica reaction-induced damage: Influence
  of alkali-silica gel properties.
\newblock \emph{Cement Concrete Res.}, 109:\penalty0 90--102, 2018.
\newblock ISSN 0008-8846.
\newblock \doi{10.1016/j.cemconres.2018.03.020}.

\bibitem[Curtis and Que(2015)]{Curtis2015}
F.~E. Curtis and X.~Que.
\newblock A quasi-{Newton} algorithm for nonconvex, nonsmooth optimization with
  global convergence guarantees.
\newblock \emph{Mathematical Programming Computation}, 7\penalty0 (4):\penalty0
  399--428, Dec. 2015.
\newblock ISSN 1867-2957.
\newblock URL \url{https://doi.org/10.1007/s12532-015-0086-2}.

\bibitem[Dai et~al.(2004)Dai, Liao, and Li]{Dai2004}
Y.-H. Dai, L.-Z. Liao, and D.~Li.
\newblock On restart procedures for the conjugate gradient method.
\newblock \emph{Numerical Algorithms}, 35\penalty0 (2):\penalty0 249--260,
  2004.

\bibitem[{de Geus} et~al.(2017){de Geus}, Vond\v{r}ejc, Zeman, Peerlings, and
  Geers]{deGeus2017}
T.~W.~J. {de Geus}, J.~Vond\v{r}ejc, J.~Zeman, R.~H.~J. Peerlings, and
  M.~Geers.
\newblock Finite strain {FFT}-based non-linear solvers made simple.
\newblock \emph{Comput. Method. Appl. M.}, 318:\penalty0 412--430, 2017.
\newblock ISSN 0045-7825.
\newblock \doi{10.1016/j.cma.2016.12.032}.

\bibitem[DeJong et~al.(2008)DeJong, Hendriks, and Rots]{Dejong2008}
M.~J. DeJong, M.~A. Hendriks, and J.~G. Rots.
\newblock Sequentially linear analysis of fracture under non-proportional
  loading.
\newblock \emph{Eng. Fract. Mech.}, 75\penalty0 (18):\penalty0 5042--5056,
  2008.

\bibitem[Eisenlohr et~al.(2013)Eisenlohr, Diehl, Lebensohn, and
  Roters]{Eisenlohr2013}
P.~Eisenlohr, M.~Diehl, R.~Lebensohn, and F.~Roters.
\newblock A spectral method solution to crystal elasto-viscoplasticity at
  finite strains.
\newblock \emph{Int. J. Plast.}, 46:\penalty0 37 -- 53, 2013.
\newblock ISSN 0749-6419.
\newblock \doi{https://doi.org/10.1016/j.ijplas.2012.09.012}.
\newblock URL
  \url{http://www.sciencedirect.com/science/article/pii/S0749641912001428}.
\newblock Microstructure-based Models of Plastic Deformation.

\bibitem[Eshelby(1957)]{Eshelby1957}
J.~D. Eshelby.
\newblock The determination of the elastic field of an ellipsoidal inclusion,
  and related problems.
\newblock \emph{P. Roy. Soc. A-Math. Phys.}, 241\penalty0 (1226):\penalty0
  376--396, 1957.
\newblock \doi{10.1098/rspa.1957.0133}.

\bibitem[Eshelby(1959)]{Eshelby1959}
J.~D. Eshelby.
\newblock The elastic field outside an ellipsoidal inclusion.
\newblock \emph{P. Roy. Soc. A-Math. Phy.}, 252\penalty0 (1271):\penalty0
  561--569, 1959.
\newblock \doi{10.1098/rspa.1959.0173}.

\bibitem[Frigo and Johnson(2005)]{Frigo2005a}
M.~Frigo and S.~G. Johnson.
\newblock The design and implementation of {FFTW3}.
\newblock \emph{Proc. IEEE}, 93\penalty0 (2):\penalty0 216--231, 2005.

\bibitem[Gallyamov et~al.(2020)Gallyamov, Ramos, Corrado, Rezakhani, and
  Molinari]{Gallyamov2020}
E.~R. Gallyamov, A.~C. Ramos, M.~Corrado, R.~Rezakhani, and J.-F. Molinari.
\newblock Multi-scale modelling of concrete structures affected by
  alkali-silica reaction: Coupling the mesoscopic damage evolution and the
  macroscopic concrete deterioration.
\newblock \emph{Int. J. Solids Struct.}, 207:\penalty0 262--278, 2020.

\bibitem[Geers et~al.(2010)Geers, Kouznetsova, and Brekelmans]{Geers2010}
M.~G.~D. Geers, V.~G. Kouznetsova, and W.~A.~M. Brekelmans.
\newblock Multi-scale computational homogenization: Trends and challenges.
\newblock \emph{J. Comput. Appl. Math.}, 234\penalty0 (7):\penalty0 2175--2182,
  2010.
\newblock ISSN 0377-0427.
\newblock \doi{10.1016/j.cam.2009.08.077}.

\bibitem[Gelb and Gottlieb(2007)]{Gelb2007}
A.~Gelb and S.~Gottlieb.
\newblock The resolution of the {Gibbs} phenomenon for {Fourier} spectral
  methods.
\newblock \emph{Advances in The {Gibbs} Phenomenon. Sampling Publishing,
  Potsdam, New York}, 2007.

\bibitem[G{\'e}l{\'e}bart and Mondon-Cancel(2013)]{Gelebart2013}
L.~G{\'e}l{\'e}bart and R.~Mondon-Cancel.
\newblock Non-linear extension of {FFT}-based methods accelerated by conjugate
  gradients to evaluate the mechanical behavior of composite materials.
\newblock \emph{Computational Materials Science}, 77:\penalty0 430--439, 2013.

\bibitem[Gottlieb and Shu(1996)]{Gottlieb1996}
D.~Gottlieb and C.-W. Shu.
\newblock On the {Gibbs} phenomenon {III}: recovering exponential accuracy in a
  sub-interval from a spectral partial sum of a piecewise analytic function.
\newblock \emph{SIAM J. Numer. Anal.}, 33\penalty0 (1):\penalty0 280--290,
  1996.

\bibitem[Gottlieb and Shu(1997)]{Gottlieb1997}
D.~Gottlieb and C.-W. Shu.
\newblock On the {Gibbs} phenomenon and its resolution.
\newblock \emph{SIAM Rev.}, 39\penalty0 (4):\penalty0 644--668, 1997.
\newblock \doi{10.1137/S0036144596301390}.

\bibitem[Hewitt and Hewitt(1979)]{Hewitt1979}
E.~Hewitt and R.~E. Hewitt.
\newblock The {Gibbs}-{Wilbraham} phenomenon: {An} episode in {Fourier}
  analysis.
\newblock \emph{Arch. Hist. Exact Sci.}, pages 129--160, 1979.

\bibitem[Hill(1963)]{Hill1963}
R.~Hill.
\newblock Elastic properties of reinforced solids: some theoretical principles.
\newblock \emph{J. Mech. Phys. Solids}, 11\penalty0 (5):\penalty0 357--372,
  1963.

\bibitem[Hill(1985)]{Hill1985}
R.~Hill.
\newblock On the micro-to-macro transition in constitutive analyses of
  elastoplastic response at finite strain.
\newblock In \emph{Mathematical proceedings of the Cambridge philosophical
  society}, volume~98, pages 579--590. Cambridge University Press, Cambridge
  University Press ({CUP}), nov 1985.
\newblock \doi{10.1017/S0305004100063787}.

\bibitem[Hobbs(1988)]{Hobbs1988}
D.~W. Hobbs.
\newblock \emph{Alkali-silica reaction in concrete}.
\newblock Thomas Telford Publishing, 1988.

\bibitem[Hsia et~al.(2017)Hsia, Zhu, and Lin]{Hsia2017}
C.-Y. Hsia, Y.~Zhu, and C.-J. Lin.
\newblock A study on trust region update rules in {Newton} methods for
  large-scale linear classification.
\newblock In \emph{Asian conference on machine learning}, pages 33--48. PMLR,
  2017.

\bibitem[Junge(2022)]{muspectre}
T.~Junge.
\newblock \url{https://gitlab.com/muspectre/muspectre}, 2022.

\bibitem[Kabel et~al.(2014)Kabel, Böhlke, and Schneider]{Kabel2014}
M.~Kabel, T.~Böhlke, and M.~Schneider.
\newblock Efficient fixed point and {Newton}-krylov solvers for {FFT}-based
  homogenization of elasticity at large deformations.
\newblock \emph{Comput. Mech.}, 54\penalty0 (6):\penalty0 1497--1514, Dec.
  2014.
\newblock ISSN 1432-0924.
\newblock URL \url{https://doi.org/10.1007/s00466-014-1071-8}.

\bibitem[Kabel et~al.(2015)Kabel, Merkert, and Schneider]{Kabel2015}
M.~Kabel, D.~Merkert, and M.~Schneider.
\newblock Use of composite voxels in {FFT}-based homogenization.
\newblock \emph{Computer Methods in Applied Mechanics and Engineering},
  294:\penalty0 168--188, 2015.

\bibitem[Khorrami et~al.(2020)Khorrami, Mianroodi, Shanthraj, and
  Svendsen]{Khorrami2020}
M.~Khorrami, J.~R. Mianroodi, P.~Shanthraj, and B.~Svendsen.
\newblock Development and comparison of spectral algorithms for numerical
  modeling of the quasi-static mechanical behavior of inhomogeneous materials.
\newblock \emph{arXiv:2009.03762}, 2020.

\bibitem[Kochmann et~al.(2018)Kochmann, Ehle, Wulfinghoff, Mayer, Svendsen, and
  Reese]{Kochmann2018}
J.~Kochmann, L.~Ehle, S.~Wulfinghoff, J.~Mayer, B.~Svendsen, and S.~Reese.
\newblock Efficient multiscale {FE}-{FFT}-based modeling and simulation of
  macroscopic deformation processes with non-linear heterogeneous
  microstructures.
\newblock In \emph{Multiscale Modeling of Heterogeneous Structures}, pages
  129--146. Springer, 2018.

\bibitem[Ladeck{\'{y}} et~al.(2021)Ladeck{\'{y}}, Pultarov{\'{a}}, and
  Zeman]{Ladecky2021a}
M.~Ladeck{\'{y}}, I.~Pultarov{\'{a}}, and J.~Zeman.
\newblock Guaranteed two-sided bounds on all eigenvalues of preconditioned
  diffusion and elasticity problems solved by the finite element method.
\newblock \emph{Appl. Math.}, 66\penalty0 (1):\penalty0 21--42, jan 2021.
\newblock \doi{10.21136/AM.2020.0217-19}.

\bibitem[Ladecký et~al.(2022)Ladecký, Leute, Falsafi, Pultarová, Pastewka,
  Junge, and Zeman]{Ladecky2022}
M.~Ladecký, R.~J. Leute, A.~Falsafi, I.~Pultarová, L.~Pastewka, T.~Junge, and
  J.~Zeman.
\newblock Optimal {FFT}-accelerated finite element solver for homogenization,
  2022.
\newblock URL \url{https://arxiv.org/abs/2203.02962}.

\bibitem[Leuschner and Fritzen(2018)]{Leuschner2018}
M.~Leuschner and F.~Fritzen.
\newblock {Fourier}-accelerated nodal solvers (fans) for homogenization
  problems.
\newblock \emph{Comput. Mech.}, 62\penalty0 (3):\penalty0 359--392, Sep 2018.
\newblock ISSN 1432-0924.
\newblock \doi{10.1007/s00466-017-1501-5}.
\newblock URL \url{10.1007/s00466-017-1501-5}.

\bibitem[Leute et~al.(2022)Leute, Ladeck{\`y}, Falsafi, J{\"o}dicke,
  Pultarov{\'a}, Zeman, Junge, and Pastewka]{Leute2022}
R.~J. Leute, M.~Ladeck{\`y}, A.~Falsafi, I.~J{\"o}dicke, I.~Pultarov{\'a},
  J.~Zeman, T.~Junge, and L.~Pastewka.
\newblock Elimination of ringing artifacts by finite-element projection in
  {FFT}-based homogenization.
\newblock \emph{J. Comput. Phys.}, page 110931, 2022.

\bibitem[Li.(2017)]{Ai2017}
A.~Li.
\newblock \emph{"Micro-architectured Metamaterials: Design and Analysis"}.
\newblock Mechanical engineering research theses and dissertations. 2.,
  Southern Methodist University, 2017.
\newblock URL \url{https://scholar.smu.edu/engineering_mechanical_etds/2}.

\bibitem[Liu et~al.(2016)Liu, Sharma, Newell, Bishop, Spear, Lester,
  Castelluccio, Bonney, and Brake]{Liu2016}
C.~Liu, P.~Sharma, P.~Newell, J.~E. Bishop, A.~Spear, B.~T. Lester, G.~M.
  Castelluccio, M.~Bonney, and M.~R. Brake.
\newblock Emergent homogonization techniques and effective dynamical
  properties.
\newblock Technical report, Sandia National Lab.(SNL-NM), Albuquerque, NM
  (United States), 2016.

\bibitem[Lucarini and Segurado(2019)]{Lucarini2019}
S.~Lucarini and J.~Segurado.
\newblock An algorithm for stress and mixed control in galerkin-based {FFT}
  homogenization.
\newblock \emph{Int. J. Numer. Methods Eng.}, 119\penalty0 (8):\penalty0
  797--805, 2019.

\bibitem[Ma et~al.(2021{\natexlab{a}})Ma, Shakoor, Vasiukov, Lomov, and
  Park]{Ma2021}
X.~Ma, M.~Shakoor, D.~Vasiukov, S.~V. Lomov, and C.~H. Park.
\newblock Numerical artifacts of fast {Fourier} transform solvers for elastic
  problems of multi-phase materials: their causes and reduction methods.
\newblock \emph{Comput. Mech.}, pages in press,~, 2021{\natexlab{a}}.
\newblock \doi{10.1007/s00466-021-02013-5}.

\bibitem[Ma et~al.(2021{\natexlab{b}})Ma, Shakoor, Vasiukov, Lomov, and
  Park]{Ma2021a}
X.~Ma, M.~Shakoor, D.~Vasiukov, S.~V. Lomov, and C.~H. Park.
\newblock Numerical artifacts of fast {Fourier} transform solvers for elastic
  problems of multi-phase materials: their causes and reduction methods.
\newblock \emph{Comput. Mech.}, 67\penalty0 (6):\penalty0 1661--1683,
  2021{\natexlab{b}}.

\bibitem[Marvi-Mashhadi et~al.(2020)Marvi-Mashhadi, Lopes, and
  LLorca]{Marvi-Mashhadi2020}
M.~Marvi-Mashhadi, C.~Lopes, and J.~LLorca.
\newblock High fidelity simulation of the mechanical behavior of closed-cell
  polyurethane foams.
\newblock \emph{J. Mech. Phys. Solids}, 135:\penalty0 103814, 2020.
\newblock ISSN 0022-5096.
\newblock \doi{https://doi.org/10.1016/j.jmps.2019.103814}.
\newblock URL
  \url{https://www.sciencedirect.com/science/article/pii/S0022509619306477}.

\bibitem[Matou{\v{s}} et~al.(2017)Matou{\v{s}}, Geers, Kouznetsova, and
  Gillman]{Matous2017}
K.~Matou{\v{s}}, M.~G. Geers, V.~G. Kouznetsova, and A.~Gillman.
\newblock {A review of predictive nonlinear theories for multiscale modeling of
  heterogeneous materials}.
\newblock \emph{J. Comput. Phys.}, 330:\penalty0 192--220, feb 2017.
\newblock ISSN 10902716.
\newblock \doi{10.1016/j.jcp.2016.10.070}.

\bibitem[Mazars and Pijaudier-Cabot(1989)]{Mazars1989}
J.~Mazars and G.~Pijaudier-Cabot.
\newblock Continuum damage theory—application to concrete.
\newblock \emph{J. Eng. Mech.}, 115\penalty0 (2):\penalty0 345--365, 1989.

\bibitem[Meng et~al.(2012)Meng, Heltsley, and Pollard]{Meng2012}
C.~Meng, W.~Heltsley, and D.~D. Pollard.
\newblock Evaluation of the {Eshelby} solution for the ellipsoidal inclusion
  and heterogeneity.
\newblock \emph{Comput. Geosci.}, 40:\penalty0 40--48, 2012.
\newblock ISSN 0098-3004.
\newblock \doi{10.1016/j.cageo.2011.07.008}.

\bibitem[Michel et~al.(2001)Michel, Moulinec, and Suquet]{Michel2001}
J.~C. Michel, H.~Moulinec, and P.~Suquet.
\newblock A computational scheme for linear and non-linear composites with
  arbitrary phase contrast.
\newblock \emph{Int. J. Numer. Methods Eng.}, 52\penalty0 (1‐2):\penalty0
  139--160, 2001.
\newblock \doi{10.1002/nme.275}.
\newblock URL \url{https://onlinelibrary.wiley.com/doi/abs/10.1002/nme.275}.

\bibitem[Milton(1995)]{Milton1995}
G.~W. Milton.
\newblock The theory of composites.
\newblock \emph{Materials and Technology}, 117:\penalty0 483--93, 1995.

\bibitem[Milton and Sawicki(2003)]{Milton2003}
G.~W. Milton and A.~Sawicki.
\newblock Theory of composites. cambridge monographs on applied and
  computational mathematics.
\newblock \emph{Appl. Mech. Rev.}, 56\penalty0 (2):\penalty0 B27--B28, 2003.

\bibitem[Mishra et~al.(2016)Mishra, Vondřejc, and Zeman]{Mishra2016}
N.~Mishra, J.~Vondřejc, and J.~Zeman.
\newblock A comparative study on low-memory iterative solvers for {FFT}-based
  homogenization of periodic media.
\newblock \emph{J. Comput. Phys.}, 321:\penalty0 151--168, 2016.
\newblock ISSN 0021-9991.
\newblock \doi{https://doi.org/10.1016/j.jcp.2016.05.041}.
\newblock URL
  \url{https://www.sciencedirect.com/science/article/pii/S0021999116301863}.

\bibitem[Mori and Tanaka(1973)]{Mori1973}
T.~Mori and K.~Tanaka.
\newblock Average stress in matrix and average elastic energy of materials with
  misfitting inclusions.
\newblock \emph{Acta Metall.}, 21\penalty0 (5):\penalty0 571--574, 1973.

\bibitem[Moshfegh and Vouvakis(2020)]{Moshfegh2020}
J.~Moshfegh and M.~N. Vouvakis.
\newblock Direct solution of {FEM} models: Are sparse direct solvers the best
  strategy?, 2020.

\bibitem[MotahariTabari and Shooshpasha(2018)]{MotahariTabari2018}
S.~MotahariTabari and I.~Shooshpasha.
\newblock Evaluation of coarse-grained mechanical properties using small direct
  shear test.
\newblock \emph{Int. J. Geotech. Eng.}, 15\penalty0 (6):\penalty0 667--679, aug
  2018.
\newblock \doi{10.1080/19386362.2018.1505310}.

\bibitem[Moulinec and Suquet(1994)]{Moulinec1994}
H.~Moulinec and P.~Suquet.
\newblock A fast numerical method for computing the linear and nonlinear
  properties of composites.
\newblock \emph{C. R. Acad. Sci. II B-Mec.}, 318:\penalty0 1417--1423, 01 1994.

\bibitem[Mura(1982)]{Mura1982}
T.~Mura.
\newblock \emph{Micromechanics of Defects in Solids}.
\newblock Kluwer Academic Publishers Group, 1982.
\newblock ISBN 978-94-011-8548-6.
\newblock \doi{10.1007/978-94-011-9306-1}.

\bibitem[Nemat-Nasser and Hori(2013)]{Nemat-Nasser2013}
S.~Nemat-Nasser and M.~Hori.
\newblock \emph{Micromechanics: overall properties of heterogeneous materials}.
\newblock Elsevier, 2013.

\bibitem[Nocedal and Wright(2006)]{Nocedal2006}
J.~Nocedal and S.~Wright.
\newblock \emph{Numerical optimization}.
\newblock Springer Science \& Business Media, 2006.

\bibitem[Norris(1985)]{Norris1985}
A.~Norris.
\newblock A differential scheme for the effective moduli of composites.
\newblock \emph{Mech. Mater.}, 4\penalty0 (1):\penalty0 1--16, 1985.

\bibitem[Pari et~al.(2022)Pari, Rots, and Hendriks]{Pari2022}
M.~Pari, J.~G. Rots, and M.~Hendriks.
\newblock Recent advancements in sequentially linear analysis ({SLA}) type
  solution procedures.
\newblock In \emph{Computational Modelling of Concrete and Concrete
  Structures}, pages 432--442. CRC Press, 2022.

\bibitem[Pijaudier-Cabot and Ba{\v{z}}ant(1987)]{Pijaudier-Cabot1987}
G.~Pijaudier-Cabot and Z.~P. Ba{\v{z}}ant.
\newblock Nonlocal damage theory.
\newblock \emph{J. Eng. Mech.}, 113\penalty0 (10):\penalty0 1512--1533, 1987.

\bibitem[Pippig(2013)]{Pippig2013}
M.~Pippig.
\newblock {PFFT}: An extension of {FFTW} to massively parallel architectures.
\newblock \emph{{SIAM} Journal on Scientific Computing}, 35\penalty0
  (3):\penalty0 C213--C236, jan 2013.
\newblock \doi{https://doi.org/10.1137/120885887}.

\bibitem[Pivovarov et~al.(2018)Pivovarov, Steinmann, and
  Willner]{Pivovarov2018}
D.~Pivovarov, P.~Steinmann, and K.~Willner.
\newblock Two reduction methods for stochastic {FEM} based homogenization using
  global basis functions.
\newblock \emph{Computer Methods in Applied Mechanics and Engineering},
  332:\penalty0 488--519, 2018.

\bibitem[Powell(1977)]{Powell1977}
M.~J.~D. Powell.
\newblock Restart procedures for the conjugate gradient method.
\newblock \emph{Math. Program.}, 12\penalty0 (1):\penalty0 241--254, 1977.

\bibitem[Prakash and Lebensohn(2009)]{Prakash2009}
A.~Prakash and R.~Lebensohn.
\newblock Simulation of micromechanical behavior of polycrystals: finite
  elements versus fast {Fourier} transforms.
\newblock \emph{Modell. Simul. Mater. Sci. Eng.}, 17\penalty0 (6):\penalty0
  064010, 2009.

\bibitem[Pultarov{\'a} and Ladeck{\'y}(2021)]{Pultarova2021}
I.~Pultarov{\'a} and M.~Ladeck{\'y}.
\newblock Two-sided guaranteed bounds to individual eigenvalues of
  preconditioned finite element and finite difference problems.
\newblock \emph{Numerical Linear Algebra with Applications}, 28\penalty0
  (5):\penalty0 e2382, 2021.

\bibitem[Roters et~al.(2013)Roters, Diehl, Shanthraj, Lebensohn, and
  Eisenlohr]{Roters2013}
F.~Roters, M.~Diehl, P.~Shanthraj, R.~Lebensohn, and P.~Eisenlohr.
\newblock A spectral method solution to crystal elastoviscoplasticity at finite
  strains.
\newblock In \emph{Plasticity’13, The 19th International Symposium on
  Plasticity \& Its Current Applications}, 2013.

\bibitem[Rots(2001)]{Rots2001}
J.~Rots.
\newblock Sequentially linear continuum model for concrete fracture.
\newblock \emph{Fracture mechanics of concrete structures}, 2:\penalty0
  831--840, 2001.

\bibitem[Rots(1988)]{Rots1988}
J.~G. Rots.
\newblock \emph{Computational modeling of concrete fracture}.
\newblock PhD thesis, Czech University of Life Sciences Prague, 1988.

\bibitem[Rots and Invernizzi(2004)]{Rots2004}
J.~G. Rots and S.~Invernizzi.
\newblock Regularized sequentially linear saw-tooth softening model.
\newblock \emph{Int. J. Numer. Anal. Methods Geomech.}, 28\penalty0
  (7-8):\penalty0 821--856, 2004.

\bibitem[Rots et~al.(2008)Rots, Belletti, and Invernizzi]{Rots2008}
J.~G. Rots, B.~Belletti, and S.~Invernizzi.
\newblock Robust modeling of rc structures with an “event-by-event”
  strategy.
\newblock \emph{Eng. Fract. Mech.}, 75\penalty0 (3-4):\penalty0 590--614, 2008.

\bibitem[Schneider(2021)]{Schneider2021}
M.~Schneider.
\newblock A review of nonlinear {FFT}-based computational homogenization
  methods.
\newblock \emph{Acta Mech.}, pages in press,~, 2021.
\newblock \doi{10.1007/s00707-021-02962-1}.

\bibitem[Schneider et~al.(2016)Schneider, Ospald, and Kabel]{Schneider2016}
M.~Schneider, F.~Ospald, and M.~Kabel.
\newblock Computational homogenization of elasticity on a staggered grid.
\newblock \emph{Int. J. Numer. Meth. Eng.}, 105\penalty0 (9):\penalty0
  693--720, 2016.
\newblock \doi{10.1002/nme.5008}.

\bibitem[Schroder(2014)]{Schroeder2014}
J.~Schroder.
\newblock \emph{A numerical two-scale homogenization scheme: the FE2-method},
  pages 1--64.
\newblock Springer Vienna, Vienna, 2014.
\newblock ISBN 978-3-7091-1625-8.
\newblock \doi{10.1007/978-3-7091-1625-8_1}.
\newblock URL \url{https://doi.org/10.1007/978-3-7091-1625-8_1}.

\bibitem[Sellier et~al.(2017)Sellier, Grimal, Multon, and
  Bourdarot]{Sellier2017}
A.~Sellier, {\'E}.~Grimal, S.~Multon, and E.~Bourdarot.
\newblock \emph{Swelling concrete in dams and hydraulic structures: DSC 2017}.
\newblock John Wiley \& Sons, 2017.

\bibitem[Steihaug(1983)]{Steihaug1983}
T.~Steihaug.
\newblock The conjugate gradient method and trust regions in large scale
  optimization.
\newblock \emph{SIAM J. Numer. Anal.}, 20\penalty0 (3):\penalty0 626--637,
  1983.

\bibitem[Swamy(1991)]{Swamy1991}
R.~Swamy.
\newblock \emph{The Alcali-Silica Raection in Concrete}.
\newblock CRC Press, 1991.

\bibitem[Vinogradov and Milton(2008)]{Vinogradov2008}
V.~Vinogradov and G.~W. Milton.
\newblock An accelerated {FFT} algorithm for thermoelastic and non-linear
  composites.
\newblock \emph{Int. J. Numer. Methods Eng.}, 76\penalty0 (11):\penalty0
  1678--1695, 2008.
\newblock \doi{10.1002/nme.2375}.
\newblock URL \url{https://onlinelibrary.wiley.com/doi/abs/10.1002/nme.2375}.

\bibitem[Vond\v{r}ejc et~al.(2014)Vond\v{r}ejc, Zeman, and Marek]{Vondrejc2014}
J.~Vond\v{r}ejc, J.~Zeman, and I.~Marek.
\newblock An {FFT}-based {G}alerkin method for homogenization of periodic
  media.
\newblock \emph{Comput. Math. Appl.}, 68\penalty0 (3):\penalty0 156--173, 2014.
\newblock ISSN 0898-1221.
\newblock \doi{10.1016/j.camwa.2014.05.014}.

\bibitem[Wicht et~al.(2019)Wicht, Schneider, and Böhlke]{Wicht2019}
D.~Wicht, M.~Schneider, and T.~Böhlke.
\newblock On quasi-newton methods in fast fourier transform-based
  micromechanics.
\newblock \emph{International Journal for Numerical Methods in Engineering},
  121\penalty0 (8):\penalty0 1665--1694, dec 2019.
\newblock \doi{10.1002/nme.6283}.

\bibitem[Willot(2015)]{Willot2015}
F.~Willot.
\newblock {Fourier}-based schemes for computing the mechanical response of
  composites with accurate local fields.
\newblock \emph{C. R. Mécanique}, 343\penalty0 (3):\penalty0 232--245, 2015.
\newblock ISSN 1631-0721.
\newblock \doi{10.1016/j.crme.2014.12.005}.

\bibitem[Yuan(2000)]{Yuan2000}
Y.-x. Yuan.
\newblock A review of trust region algorithms for optimization.
\newblock In \emph{Iciam}, volume~99, pages 271--282, 2000.

\bibitem[Zeman et~al.(2010)Zeman, Vond{\v{r}}ejc, Nov{\'a}k, and
  Marek]{Zeman2010}
J.~Zeman, J.~Vond{\v{r}}ejc, J.~Nov{\'a}k, and I.~Marek.
\newblock Accelerating a {FFT}-based solver for numerical homogenization of
  periodic media by conjugate gradients.
\newblock \emph{J. Comput. Phys.}, 229\penalty0 (21):\penalty0 8065--8071,
  2010.

\bibitem[Zeman et~al.(2017)Zeman, de~Geus, Vond\v{r}ejc, Peerlings, and
  Geers]{Zeman2017}
J.~Zeman, T.~W.~J. de~Geus, J.~Vond\v{r}ejc, R.~H.~J. Peerlings, and M.~G.~D.
  Geers.
\newblock A finite element perspective on nonlinear {FFT}-based micromechanical
  simulations.
\newblock \emph{Int. J. Numer. Meth. Eng.}, 111\penalty0 (10):\penalty0
  903--926, 2017.
\newblock \doi{10.1002/nme.5481}.

\end{thebibliography}

\end{document}